\documentclass[11pt]{article}

\usepackage{amsmath}
\usepackage{amssymb}
\usepackage{latexsym}
\usepackage{epsfig}

\def\F{\mathbb{F}}

{\catcode `\@=11 \global\let\AddToReset=\@addtoreset}
\AddToReset{equation}{section}

\usepackage{graphicx}
\usepackage{color}

\usepackage[latin1]{inputenc}

\def\I{\mathbb{I}}

\newtheorem{theorem}{Theorem}
\newtheorem{lemma}{Lemma}
\newtheorem{corollary}{Corollary}
\newtheorem{proposition}{Proposition}
\newtheorem{@definition}{Definition}
\newenvironment{definition}{\begin{@definition}\rm}{\end{@definition}}
\newtheorem{@remark}{Remark}
\newenvironment{remark}{\begin{@remark}\rm}{\end{@remark}}
\newtheorem{@example}{Example}
\newenvironment{example}{\begin{@example}\rm}{\end{@example}}

\newcommand{\beqn}{\begin{displaymath}}
\newcommand{\eeqn}{\end{displaymath}}
\newcommand{\beq}{\begin{equation}}  
\newcommand{\eeq}{\end{equation}}

\def\mathsf{\bf}
\def\N{\mathbb{N}}

\def\R{\mathbb{R}}
\def\C{\mathbb{C}}

\def\I{\mathbb{I}}
\def\i{\mathrm i}
\def\d{\mathrm d}
\def\e{\mathrm e}

\def\E{\mathrm E}
\def\P{\mathrm P}

\def\text{\mbox}

\def\1{{\bf 1}}

\newcommand{\Leb}{{\rm Leb}}

\newcommand{\Var}{\mbox{Var}\,}
\newcommand{\Cov}{\mbox{Cov}\,}

\newcommand{\nn}{\nonumber}
\newcommand{\noi}{\noindent}

\newcommand{\mbt}{\boldsymbol{t}}

\newcommand{\mbs}{\boldsymbol{s}}

\newcommand{\mbu}{\boldsymbol{u}}

\newcommand{\mbx}{\boldsymbol{x}}

\newcommand{\mbz}{\boldsymbol{z}}

\newcommand{\mbgamma}{\boldsymbol {\gamma}}

\definecolor{vp}{rgb}{0.55, 0.71, 0.0}

\def\eq2{
\stackrel{\small \rm mod \,2}{=}}

\def\n2{
\stackrel{\small \rm mod \,2}{\neq}}

\def\limd{\renewcommand{\arraystretch}{0.5}
\begin{array}[t]{c}
\stackrel{\rm d}{\longrightarrow} \\
\end{array}\renewcommand{\arraystretch}{1}}

\def\limfdd{\renewcommand{\arraystretch}{0.5}
\begin{array}[t]{c}
\stackrel{\rm fdd}{\longrightarrow} \\
\end{array}\renewcommand{\arraystretch}{1}}

\def\eqd{\renewcommand{\arraystretch}{0.5}
\begin{array}[t]{c}
\stackrel{\rm d}{=} \\
\end{array}\renewcommand{\arraystretch}{1}}

\def\1{{\bf 1}}
\def\0{{\bf 0}}

\def\neqd{
	\begin{array}[t]{c}
		\stackrel{{\rm d}}{\neq}
\end{array}}

\def\Cov{\mathrm{Cov}\,}

\begin{document}

\title{Limit distribution of the sample volume fraction of  Boolean set}
\author{Hermine Bierm\'e$^1$, \  Olivier Durieu$^1$ \ and \
Donatas Surgailis$^2$}
\date{\today
\\  \small
\vskip.2cm
$^1$ Universit\'e de Tours,
Universit\'e d'Orl\'eans, CNRS, IDP UMR 7013, Tours, France\\
$^2$Vilnius University, Faculty of Mathematics and Informatics,
Naugarduko 24, 03225  Vilnius, Lithuania \\
}
\maketitle

\begin{abstract}

\medskip

We study the limit distribution of the volume fraction estimator $\widehat p_{\lambda, A}$
defined as the Lebesgue  measure of the intersection
$\mathcal{X}\cap (\lambda A)$  of a random set $\mathcal{X}$ with a large observation set $\lambda A$, divided by the Lebesgue  measure of $\lambda A$,
as $\lambda \to \infty$,
for a Boolean set  $\mathcal{X}$ formed by uniformly scattered random grains $\Xi  \subset \mathbb{R}^\nu$. We obtain general conditions
on the generic grain set $\Xi$ under which $\widehat p_{\lambda, A}$  has an $\alpha$-stable limit distribution with index
$1 < \alpha \le  2$. A large class of Boolean models with randomly homothetic grains
satisfying these conditions
is introduced. We also discuss the limit distribution of the sample volume  fraction of a Boolean set
 observed
on a large subset of a  $\nu_0$-dimensional hyperplane of  $\mathbb{R}^\nu$ ($1 \le \nu_0 \le \nu -1$).
Similar results are also obtained for more general excursion sets of Boolean models.

\end{abstract}

{\small

\noi {\it Keywords:} Poisson process, Boolean set,
random grain model, long-range dependence, volume fraction
estimator on hyperplane,
stable limit distribution, randomly homothetic grain, excursion set.

\noi {\it 2020 Mathematics Subject Classification:} Primary, 60F05, 60G60, 60D05; secondary, 60G55.
}

\section{Introduction}

Volume fraction estimation, also called area fraction in 2D or porosity in porous media, is a fundamental problem in stochastic geometry, especially when working with Boolean models \cite{ baccelli2024random, sto1989}. A Boolean model is obtained by considering
 $\{\mbu_j; j \ge 1 \} $  a stationary Poisson process on $\R^\nu$ with unit intensity (called germ process) and
$\{\Xi, \Xi_j; j \ge 1 \} $  an independent and identically distributed (i.i.d.) sequence of random sets  in $ \R^\nu$ (called grains), independent of $\{\mbu_j; j \ge 1 \} $.
A {\it Boolean set} is defined as the union of all the grains:
\begin{eqnarray} \label{bool}
{\mathcal X} := \bigcup_{j=1}^\infty (\mbu_j + \Xi_j) \ \subset \ \R^\nu.
\end{eqnarray}
The Boolean model is the most important coverage model in stereology and stochastic geometry, see  \cite{sto1989}.
Rigorous definitions of  \eqref{bool} and random sets are given later.

A closely related  {\it random grain model}   is defined  as the superposition of the  indicator functions of grains $\mbu_j + \Xi_j$, viz.,
\begin{eqnarray}\label{rG1}
X(\mbt) = \sum_{j=1}^\infty \I (\mbt \in (\mbu_j + \Xi_j)), \quad \mbt \in  \R^\nu.
\end{eqnarray}
Considering randomly dilated balls for grains yields {\it random balls models} that have been considered in \cite{bier2010} and generalized in \cite{bier2018,breton2009rescaled}.
The Boolean  set in \eqref{bool} can be identified with its indicator function $\widehat X(\mbt) := \I (\mbt \in {\mathcal X})$,  $\mbt \in \R^\nu $, which is
a simple nonlinear transformation of the linear random field in
\eqref{rG1}
\begin{eqnarray}\label{bool3}
\widehat X(\mbt) = X(\mbt) \wedge 1,
\end{eqnarray}
where $a\wedge b=\min(a,b)$ for real numbers $a, b$.
The basic assumption guaranteeing  the finiteness of \eqref{rG1} is
\begin{eqnarray} \label{Wgraincond}
\mu := \E \, \Leb_\nu (\Xi) < \infty.
\end{eqnarray}
In this paper, a random closed set satisfying \eqref{Wgraincond} is called a random grain and the random grain model in \eqref{rG1} is well-defined.
 It has marginal Poisson
 distributions with mean $\mu $ and  non-negative covariance function
\begin{eqnarray} \label{covRG}
{\rm Cov}(X (\boldsymbol{0}), X(\mbt))
&=&
\E \Leb_\nu (\Xi \cap (\Xi - \mbt)) \ge 0, \quad \mbt \in \R^\nu.
\end{eqnarray}

The {\it volume fraction} of the stationary Boolean set $\mathcal{X}$ in \eqref{bool} is
the mean of the volume of $\mathcal{X}$ in the unit cube $(\boldsymbol{0}, \1]:=(0,1]^\nu$, that is
\begin{eqnarray}\label{p1}
p := \E \Leb_\nu (\mathcal{X} \cap (\boldsymbol{0}, \1]) = \int_{]\boldsymbol{0}, \1]} \E \widehat X(\mbt) \d \mbt = \E \widehat X(\boldsymbol{0}) =
1 - \P(X(\boldsymbol{0})=0),
\end{eqnarray}
leading to
\begin{eqnarray}\label{p2}
p = 1 - \e^{- \mu}.
\end{eqnarray}
The volume fraction is the most important parameter of a Boolean set,
the analog of the mean (expectation) of a stationary process on $\R^\nu$.
By stationarity,  $p= \E \Leb_\nu (\mathcal{X}\cap A)/ \Leb_\nu (A)$ for any Borel set $A$ with $0< \Leb_\nu (A) < \infty$.
The natural estimator of $p$ from observations of ${\mathcal X}$ on a large inflated set $\lambda A \subset \R^\nu$ is the ratio
\begin{eqnarray}\label{volumef}
\widehat p_{\lambda, A} &:=&\frac{ \widehat X_\lambda (A)}{\Leb_\nu (\lambda A)},
\end{eqnarray}
called the {\it sample volume fraction},
where $\widehat X_\lambda (A) :=  \Leb_\nu ({\mathcal X} \cap \lambda A) $ is the volume of the intersection of the Boolean set with
$\lambda A$. Then
\begin{eqnarray}\label{volumef1}
\lambda^\nu (\widehat p_{\lambda, A}-p)&=&\frac{ \widehat X_\lambda (A)-  \E \widehat X_\lambda (A)}{\Leb_\nu (A)}
\end{eqnarray}
and finding the limit distribution of the left hand side reduces to  that of the numerator of the last fraction.

A stationary random field $Y = \{Y(\mbt);  \mbt \in \R^\nu\} $ with finite variance and covariance $r_Y(\mbt) = {\rm Cov}(Y(\boldsymbol{0}), $ $ Y(\mbt))$
is said {\it long-range dependent} if $\int_{\R^\nu} |r_Y(\mbt)| \d \mbt = \infty $ and
{\it short-range dependent} if $\int_{\R^\nu} |r_Y(\mbt)| \d \mbt < \infty$, with $\int_{\R^\nu} r_Y(\mbt) \d \mbt \ne 0$ \cite{lah2016, pils2024, sur2024}.
Under \eqref{Wgraincond}, the random grain model $X$ in \eqref{rG1}
has non negative covariance given in \eqref{covRG} and we see that
\begin{eqnarray*}
\int_{\R^\nu}{\rm Cov}(X (\boldsymbol{0}), X(\mbt)) \d \mbt
&=&\E \int_{\R^\nu }
\I (\mbs  \in   \Xi) \d \mbs \int_{\R^\nu} \I(\mbs + \mbt \in \Xi) \d \mbt\\
&=&\E \Leb_\nu (\Xi)^2.
\end{eqnarray*}
Therefore, it is long-range dependent if  $\E \Leb_\nu (\Xi)^2 = \infty $ and
short-range dependent if  $\E \Leb_\nu (\Xi)^2 < \infty $.

Clearly, a random grain model is long-range dependent if the distribution function of
$\Leb_\nu(\Xi)$ is regularly varying with index $1< \alpha < 2$. In our work,
we assume the stronger condition
\begin{eqnarray}\label{1volf}
&\P( \Leb_\nu(\Xi) >  x) \ \sim \  c_{\Xi}\, x^{-\alpha}, \quad  x \to \infty, \quad  \text{for some }  c_{\Xi} >0 \text{ and }
1 < \alpha < 2.
\end{eqnarray}
Extension of our results to the more general case when the asymptotic in
\eqref{1volf} contains a slowly varying factor is an open problem.

In most of the literature on long-range dependent random grain models  \cite{bier2018, bier2010,  kaj2007, sur2024}, it is assumed that the randomness of $\Xi $ is due to a dilation of a  {\it deterministic} set  $\Xi^0 $ by a random factor $R^{1/\nu}$, viz.,
\begin{eqnarray}\label{Theta0}
\Xi = R^{1/\nu} \, \Xi^0,
\end{eqnarray}
where $R >0$ is a random variable with regularly decaying $\alpha$-tail,  $\alpha \in (1,2)$.
For $\Xi$ in \eqref{Theta0}, the results in \cite{sur2024} imply
that the sample volume fraction in \eqref{volumef} has an asymmetric $\alpha$-stable limit distribution, for arbitrary bounded
Borel set $A$. The definition \eqref{Theta0} comprises a very special class of random set with all grains homothetic to each other. The present paper
extends the result in \cite{sur2024} to a much more general grain class.

One of the major results of this work is Theorem \ref{thm1}, implying Corollary \ref{Cor:thm1}
which says that
condition \eqref{1volf} together with
\begin{eqnarray}\label{2volf}
\E \Leb_\nu\big(\Xi\cap   \{|\mbt| > \lambda\}\big)&=&o(\lambda^{(1-\alpha)\nu/\alpha}), \qquad \lambda \to \infty
\end{eqnarray}
imply  that the centered and normalized sample volume fraction
$\widehat p_{\lambda, A}$
has an $\alpha$-stable limit distribution. We note that
conditions \eqref{1volf} and \eqref{2volf}  involve the  Lebesgue measure of
$\Xi$ and $\Xi\cap   \{|\mbt| > \lambda\} $ alone and do not
impose any structural assumptions on the Boolean set in contrast to  \eqref{Theta0}; moreover, the sufficient condition \eqref{2volf} is sharp
in the sense that the exponent $(1-\alpha)\nu/\alpha$ cannot
be improved in general,
see the examples in Section~\ref{sec:disc}.
Theorem \ref{thm1} and Corollary \ref{Cor:thm1} refer to the long-range dependent random grain model; in the short-range dependent case $\E \Leb_\nu (\Xi)^2 < \infty $,  without any additional conditions,
we also prove a central limit theorem in Theorem~\ref{thm2}, recovering a classical central limit theorem for the sample volume fraction (Corollary \ref{Cor:thm2}).
Note that
Gaussian limits for estimators of $p$ were already obtained in \cite[Theorem~3.5]{Hall88}, see also \cite{badd1980,
mase1982, molch_stoy1994} under more stringent assumptions on $\Xi$ and the observation set.

We also consider excursion sets of random grain random fields given for any $k=1,2,\ldots$ and  $u\in [k-1,k)$ by
$$\{\mbt; X(\mbt)> u\}=\{\mbt; X(\mbt)\ge k\}.$$
The Boolean set ${\mathcal X}$ corresponds to the excursion set of any level $u\in [0,1)$.
 There is a growing interest in the study of the mean geometry of excursion sets of random fields in view of their links with extremal properties  \cite{AdlerReview00, AdlerTaylorBook}. Whereas most work in this areas refer to  Gaussian stationary random fields, there are also results  concerning shot noise random fields \cite{bier2020,lachieze2019normal,trapp2024geometric} or more generally, infinitely divisible random fields. Central limit theorems for the excursion set volumes under various weak dependence conditions have been investigated in \cite{Spodarev12}.
A new notion of short-range dependence based on excursion sets was also introduced in \cite{kulik2021long}, see our Remark \ref{rem:srd}.

The above facts may serve as
a motivation for extending our study to the Boolean fields corresponding to excursion sets
\begin{equation}\label{Xk}
\widehat{X}_k(\mbt)=\I(X(\mbt)\ge k), \quad k \ge 1,
\end{equation}
see Theorems \ref{thm1} and \ref{thm2}, whose proofs
use the crucial relation
in \eqref{bool3} between the Boolean set and the random grain model
and Charlier expansion of Poisson subordinated functionals
discussed in \cite{sur2024}.

The second part of this paper is devoted to
estimation of the volume fraction
from observations
on a hyperplane
\begin{eqnarray*}
H_{\nu_0}
:=
\{\mbt\in \R^\nu: \langle  \mbt, \mbgamma_i  \rangle = 0, i=1, \cdots, \nu - \nu_0 \} \subset \R^\nu
\end{eqnarray*}
of dimension $\nu_0 \in \{1, \cdots, \nu-1\} $, determined by $\nu - \nu_0$ vectors $\mbgamma_i \in \R^\nu, i=1, \cdots, \nu - \nu_0$.
This question is important in stereological applications and has been discussed in the literature for specific
Boolean sets.
The corresponding
estimator is naturally defined as
\begin{eqnarray}\label{volumefgamma}
\widehat p_{\lambda, A} (H_{\nu_0})&:=&\frac{\Leb_{\nu_0}({\mathcal X} \cap H_{\nu_0} \cap (\lambda A))}
{\Leb_{\nu_0}(H_{\nu_0} \cap (\lambda A) )}.
\end{eqnarray}
Note that $\E \widehat p_{\lambda, A} (H_{\nu_0}) =
\E \Leb_{\nu_0}({\mathcal X} \cap H_{\nu_0} \cap (\lambda A))/ \Leb_{\nu_0}(H_{\nu_0} \cap (\lambda A)) $ and
$$
\E \Leb_{\nu_0}({\mathcal X} \cap H_{\nu_0} \cap (\lambda A))
= \int_{H_{\nu_0} \cap (\lambda A)} \E (X(\mbt) \wedge 1) \d_{\nu_0} \mbt =
p   \Leb_{\nu_0}(H_{\nu_0)} \cap (\lambda A)),
$$
see \eqref{p1}, \eqref{p2}, so that
\eqref{volumefgamma} is an unbiased estimator of $p$ for any $H_{\nu_0}$ and
$\nu_0 \in \{1, \cdots, \nu-1\}$ (a surprising but simple consequence of stationarity).
Then, what can we say about the limit distribution of $\widehat p_{\lambda, A}(H_{\nu_0})$? Is it the same as in Corollary \ref{Cor:thm1} (for  long-range dependent
Boolean set), or different? If so, how does it depend on $H_{\nu_0}$ and especially, on the dimension
$\nu_0$ of this hyperplane?

The above questions involve the $\nu_0$-dimensional Lebesgue measure
of the intersections $\Xi \cap H_{\nu_0} \cap (\lambda A)$ which
may be very singular random sets in general.
It seems that further assumptions
in addition to those in Theorems \ref{thm1} and \ref{thm2}
on $\Xi$ are needed to consider the behavior on the hyperplanes. In the present work, we introduce a class
of {\it randomly  homothetic sets}
$\Xi$ having the form as in \eqref{Theta0} except that
$\Xi^0$ is a {\it random}  bounded closed set, independent of $R$
with tail behavior as in  \eqref{2volf}, $1 < \alpha < 2$.
We prove (see  Theorem~\ref{thm3} and Corollary~\ref{cor3} for precise formulations) that for such $\Xi$, after centering and normalization,
the estimator $\widehat p_{\lambda, A} (H_{\nu_0}) $  in \eqref{volumefgamma} has, after appropriate centering and normalization, an $\alpha_{0}$-stable limit distribution when $\alpha < 1 + \frac{\nu_0}{\nu}$  with
\begin{eqnarray*}
\alpha_{0} = 1 + \frac{\nu}{\nu_0}(\alpha -1) \in (\alpha,2),
\end{eqnarray*}
and a Gaussian limit distribution when  $\alpha > 1 + \frac{\nu_0}{\nu}$.

The rest of the paper is organized as follows. In Section \ref{sec:scalingLimit} we obtain
limit distributions of integrals $\widehat X_{\lambda,k} (\phi) := \int_{\R^\nu} \phi(\mbt/\lambda) \I(X(\mbt) \ge k) \d \mbt $ of random grain model in \eqref{rG1}, for any $k  = 1,2, \cdots$ and
any $\phi$  from a class $\Phi$ of test functions under
assumptions \eqref{1volf} and \eqref{2volf}, which include
the limit of the sample volume fraction in \eqref{volumef}
as the special case $k=1$, $\phi(\mbt) = \I(\mbt \in A)$.
Section \ref{sec:homothetic} introduces the randomly homothetic random grain model and discusses
its long-range dependence properties. Section \ref{sec:hyperplan} is devoted to the hyperplane observations, stating and proving
Theorem \ref{thm3}.
Finally, numerical illustrations are given in Section~\ref{sec:num}.

\smallskip

{\it Notation.}
  In what follows, $C$ denotes  a generic positive constant which may be different at different locations.
	We write $\limd$, $\eqd$, $\neqd$ for the weak convergence, equality, and inequality of distributions, $\limfdd$ for the finite dimensional convergence of distributions. We set
	$\1 := (1,\cdots, 1) \in \R^\nu$, $\0 := (0,\cdots, 0) \in \R^\nu$, and
  $\|f\|_\alpha := (\int_{\R^\nu} |f(\mbu)|^\alpha \d \mbu)^{1/\alpha}$, $\alpha >0$.
	For a Borelian set $A\subset\R^\nu$, $\I(A)$ stands for its indicator function  and $\Leb_\nu (A)$ for its Lebesgue measure.

\section{Scaling limits for excursion sets of random grain models}\label{sec:scalingLimit}

 It is usual in stochastic geometry to consider grains as closed random sets and we denote by $(\F(\R^\nu),\mathcal{B}(\F(\R^\nu)))$ the measurable space of closed subsets of $\R^\nu$, endowed with the $\sigma$-algebra induced by the Fell topology (see \cite{baccelli2024random} Chapter 9 for instance). Let $(\Omega,\mathcal A,\P)$ be a complete probability space. Assuming that $\Xi$ is a random closed set means that $\Xi:(\Omega,\mathcal A)\rightarrow (\F(\R^\nu),\mathcal{B}(\F(\R^\nu)))$ is measurable and we denote by $\P_{_{\Xi}}$ its probability distribution.

Then our random grain model $X$ given by \eqref{rG1}
   admits the Poisson integral representation
\begin{eqnarray}\label{RG1}
X(\mbt)
&=&\int_{\R^\nu \times \F(\R^\nu)}
\I (\mbt \in (\mbu + m)) {\cal N}(\d \mbu, \d  m) \\
&=&\mu + \int_{\R^\nu \times \F(\R^\nu)}
\I (\mbt \in (\mbu +  m)) \widetilde {\cal N}(\d \mbu, \d  m),
 \quad \mbt \in  \R^\nu, \nn
\end{eqnarray}
where ${\cal N}(\d \mbu, \d  m)$ is a Poisson random measure with intensity $\d \mbu \P_{_\Xi}(\d m)$, and $  \widetilde {\cal N}(\d \mbu, \d m)
=  \widetilde {\cal N}(\d \mbu, \d m) - \E  \widetilde {\cal N}(\d \mbu, \d m)$.
The random field in \eqref{RG1} is a Poisson shot noise field with kernel function given by $g_{m}=\I_{m}$ for $m\in \F(\R^\nu)$ (see \cite{baccelli2024random} Section 2.4 for instance). Under assumption \eqref{Wgraincond}, we can view our random grain model $X$ as a random variable in $L^1_{loc}(\R^\nu)$, the space of locally integrable functions, endowed with its Borel $\sigma$-algebra induced by its natural topology.

For any $k\ge 1$ we will consider the excursion set
$$\{X \ge k\}:=\{\mbt \in\R^\nu; X(\mbt)\ge k\}.$$
Note that since $X$ is a random variable with values in $L^1_{loc}(\R^\nu)$ it is also the case of $\I(X \ge k)$. It follows that $\{X \ge k\}$ is a random measurable set
as introduced in \cite{galerne2015random} (see also Section 4 of  \cite{lachieze2019normal}). Let us denote $\widehat X_k(\mbt) := \I(X(\mbt) \ge k), $ and for $\phi\in \Phi$,
\begin{eqnarray}\label{def:Xlambda}
\widehat X_{\lambda, k} (\phi) :=  \int_{\R^\nu} \phi(\mbt/\lambda)
\widehat X_k(\mbt) \d \mbt, \quad X_\lambda (\phi) :=  \int_{\R^\nu} \phi(\mbt/\lambda)
X(\mbt) \d \mbt,
\end{eqnarray}
where
\begin{equation}\label{Phi}
   \Phi := L^1(\R^\nu) \cap L^\infty (\R^\nu)
\end{equation}
ensures the a.s.\ absolute convergence of the integrals in \eqref{def:Xlambda} and the fact that both $X_\lambda (\phi)$ and  $\widehat X_{\lambda, k} (\phi)$ have finite expectation.
Note that  $\widehat X_{\lambda} (A)=  \int_{\lambda A} \widehat X(\mbt) \d \mbt = \widehat X_{\lambda, 1} (\phi)$ for $\phi=\I(A)$.
Thus, as seen from \eqref{volumef1}, the limit of $\widehat p_{\lambda, A}$ reduces to that of $\widehat X_{\lambda,1} (\phi)$ for $\phi=\I(A)$.

\subsection{Long-range dependence case}

For $\alpha\in(1,2)$, let us write
$L_\alpha(\phi) = \int_{\R^\nu} \phi(\mbt) L_\alpha (\d \mbt)$ for the $\alpha$-stable stochastic integral
with log-characteristic function
\begin{eqnarray}\label{limjRG}
j(\theta; \phi)&:=&\log  \E \e^{\i \theta L_\alpha(\phi)} \ = \
\i c_\Xi \int_{\R^\nu} \theta \phi(\mbs) \big\{\int_{\R_+}
\big(\e^{\i \theta \phi (\mbs)x} -1 \big)x^{-\alpha} \d x\big\} \d \mbs,\,\theta\in\R,
\end{eqnarray}
which is well-defined for any $\phi \in L^\alpha(\R^\nu)$,
hence also for $\phi \in \Phi$ in \eqref{Phi}.

\begin{theorem} \label{thm1}  Let $X$ be the random grain model in \eqref{rG1} with generic grain
satisfying \eqref{1volf} and \eqref{2volf} for $1 < \alpha < 2 $.
Then for any $\phi \in \Phi$ and  $X_\lambda (\phi)$, $\widehat X_{\lambda, k} (\phi)$, $k\ge 1$, given in \eqref{def:Xlambda}, one has
\begin{eqnarray}\label{limvolX}
\lambda^{-\nu/\alpha}
\big(X_{\lambda}(\phi)
- \E X_{\lambda}(\phi) \big)
&\limd&L_\alpha (\phi)
\end{eqnarray}
and
\begin{eqnarray}\label{limvolX2}
\left\{\lambda^{-\nu/\alpha}
\big(\widehat X_{\lambda,k}(\phi)
- \E \widehat X_{\lambda,k}(\phi) \big); k\ge 1\right\}
&\overset{fdd}{\longrightarrow}& \left\{\e^{-\mu} \frac{\mu^{k-1}}{(k-1)!} L_\alpha (\phi); k\ge 1\right\},
\end{eqnarray}
as $\lambda\to\infty$, where $L_\alpha(\phi)$ admits the log-characteristic function given by \eqref{limjRG}.
\end{theorem}

\begin{remark}\label{rem1}
Note that the limit in \eqref{limvolX2} corresponds to a deterministic multiple of the same $L_\alpha(\phi)$ from \eqref{limvolX}. We shall see in the proof that this follows from the fact that the leading term is the first order term in the Charlier decomposition of $\widehat X_{k}(\mbt)$.
\end{remark}

As a particular case, recalling \eqref{volumef1}, we obtain the following corollary for the volume fraction estimator of the Boolean model.

\begin{corollary}\label{Cor:thm1} Under the assumptions of Theorem \ref{thm1}, for an arbitrary bounded Borel set  $A \subset \R^\nu$, $\Leb_\nu (A) >0$,
\begin{equation*}
\lambda^{\nu -(\nu/\alpha)}(\widehat p_{\lambda, A} -
p) \limd
\e^{-\mu} L_{\alpha}(A)/\Leb_\nu(A),  \quad \mbox{ as }  \lambda \to \infty,
\end{equation*}
where $ L_{\alpha}(A) = \int_{A} L_\alpha (\d \mbt) $ has $\alpha$-stable distribution with characteristic function
\begin{eqnarray}\label{cf:cor1}
\E\e^{\i \theta L_\alpha(A)} =
\exp\Big\{ \i c_\Xi \theta \Leb_\nu(A) \int_{\R_+}
\big(\e^{\i \theta x} -1 \big)x^{-\alpha} \d x\Big\} .
\end{eqnarray}

\end{corollary}

\noi {\bf Proof of Theorem~\ref{thm1}.} The proof is accomplished in two steps. The first step is more involved and consists in proving
 the $\alpha$-stable limit in \eqref{limvolX}
using conditions \eqref{1volf}-\eqref{2volf} and the
characteristic function of the stochastic integral in \eqref{RG1}. The second step extends \eqref{limvolX} to
\eqref{limvolX2}, using the  Charlier expansion of the indicator function $\I (x \ge k)  $ as in
\cite[Corollary 1]{sur2024}.

\smallskip

\noi {\it Step 1: proof of \eqref{limvolX}.  }
Let $j_\lambda (\theta; \phi) :=   \log \E \exp \{ \i  \theta \lambda^{-\nu/\alpha}
 (X_{\lambda}(\phi) - \E X_{\lambda}(\phi))  \}$, with $\theta\in\R$, $\phi \in \Phi$,
 and $\Psi(z) := \e^{\i z} -1 - \i z$,  $z \in \R$.
Then
\begin{eqnarray}\label{jlambda}
j_\lambda(\theta; \phi)
&=&\int_{\R^\nu} \E \Psi \Big(\frac{\theta}{\lambda^{\nu/\alpha}}
\int_{\R^\nu}\phi(\mbt/\lambda) \I (\mbt - \mbs \in \Xi) \d \mbt \Big) \d \mbs \\
&=&\lambda^\nu \int_{\R^\nu} \E \Psi \Big(\frac{\theta}{\lambda^{\nu/\alpha}} \int_{\R^\nu}
\phi\big(\frac{\mbt}{\lambda} +  \mbs\big) \I (\mbt \in \Xi) \d \mbt \Big) \d \mbs. \nn
\end{eqnarray}
The intuitive argument leading to $j_\lambda (\theta; \phi)  \to j (\theta; \phi) $
uses the observation that
\begin{equation}\label{phiproblem}
\int_{\R^\nu}\phi\big(\frac{\mbt}{\lambda} +  \mbs\big) \I (\mbt \in \Xi) \d \mbt \ \to \ \phi(\mbs) \Leb_\nu(\Xi)\;\text{a.s.}
\end{equation}
as $\lambda\to\infty$, at each continuity point $\mbs$ of $\phi(\cdot)$.
Using integration by parts and the tail condition in \eqref{1volf} we see that
\begin{eqnarray}
\lambda^\nu \E \Psi \Big(\frac{\theta \phi(\mbs)}{\lambda^{\nu/\alpha}} \Leb_\nu(\Xi) \Big)
&=&\i \theta \phi(\mbs)  \int_{\R_+} (\e^{\i  \theta \phi(\mbs) x} -1) \lambda^\nu \P(\Leb_\nu(\Xi) > x \lambda^{\nu/\alpha} ) \d x  \nn \\
&\sim&\i c_\Xi \theta \phi(\mbs)
 \int_{\R_+} (\e^{\i  \theta \phi(\mbs) x} -1) x^{-\alpha} \d x. \label{lalaj}
\end{eqnarray}
Hence, if the inner  integral in \eqref{jlambda}
can be replaced by the right hand side of \eqref{phiproblem}, i.e.\
$j_\lambda(\theta; \phi)$ can be replaced by
 \begin{eqnarray*}
\tilde j_\lambda(\theta; \phi)
&:=&\lambda^\nu \int_{\R^\nu} \E \Psi \Big(\frac{\theta\phi(\mbs)}{\lambda^{\nu/\alpha}}\Leb_\nu(\Xi) \Big) \d \mbs,
\end{eqnarray*}
then, the statement of the theorem
will follow rather easily.
A rigorous  justification of the above argument using condition \eqref{2volf} is somewhat involved.
We face two difficulties. Firstly, $\phi$ is not necessarily continuous and secondly, even if
it is,  the convergence in \eqref{phiproblem} does not necessarily hold for large $|\mbt| = O(\lambda) $.

The classical Lusin's theorem states that each measurable function $\phi$ is nearly continuous. More precisely,
for any $r >0$, $\epsilon >0$, there is a measurable set
$U_{\epsilon,r} \subset B_r := \{\mbu \in \R^\nu : | \mbu| < r \} \subset \R^\nu $ such that
$\phi $ restricted to $U_{\epsilon, r}$ is continuous and
$\Leb_\nu (B_r \setminus U_{\epsilon,r}) < \epsilon $. Accordingly, denote
$$
\xi_\lambda (\mbu)
:= \lambda^{-\nu/\alpha}
\int_{\R^\nu}\phi(\mbt/\lambda + \mbu) \I (\mbt  \in \Xi) \d \mbt, \qquad j_{\lambda}(\theta; \phi, U_{\epsilon,r})
:=
\lambda^\nu \int_{U_{\epsilon,r}} \E
\Psi (\theta \xi_\lambda (\mbu) ) \d \mbu.
$$
Note that $j_\lambda(\theta; \phi) - j_\lambda (\theta; \phi, U_{\epsilon,r}) = j_\lambda(\theta; \phi, U_{\epsilon,r}^c)$ for $U_{\epsilon,r}^c := \R^\nu \setminus U_{\epsilon,r}$.
Therefore, to obtain the convergence
$j_\lambda (\theta; \phi) \to j(\theta; \phi)$ as $\lambda \to \infty$ it suffices to show that the two following relations hold:
\begin{align}\label{jnew}
&\lim_{\epsilon \to  0,\, r \to \infty} \limsup_{\lambda \to \infty} |j_\lambda(\theta; \phi,  U_{\epsilon, r}^c)
| = 0, \\  \quad
&\forall \epsilon, r >0, \qquad \lim_{\lambda\to \infty}j_\lambda(\theta; \phi, U_{\epsilon, r}) = j(\theta; \phi, U_{\epsilon, r}),
\label{jnew1}
\end{align}
where (compare with \eqref{limjRG})
$$
j(\theta; \phi, U_{\epsilon,r}) := \int_{U_{\epsilon,r}} \i c_\Xi \theta \phi(\mbu) \Big\{ \int_{\R_+} (\e^{\i \theta \phi(\mbu) x} -1) x^{-\alpha} \d x \Big\} \d \mbu.
$$

We first consider \eqref{jnew}. We recall that $\phi \in L^1(\R^\nu) \cap L^\infty (\R^\nu)$ and, for $r>0$, $\epsilon>0$, we introduce the integrals
$$
\Phi_p (\mbt)
:= \Big(\int_{U^c_{\epsilon,r}} |\phi(\mbt + \mbu)|^p \d \mbu \Big)^{1/p}, \qquad p =1,2,
$$
satisfying
$(\Phi_2 (\mbt))^2 \le \Phi_1 (\mbt)  \| \phi \|_\infty \le  \| \phi \|_1 \| \phi \|_\infty$ for each $\mbt \in \R^\nu$.
In particular, for each
$ |\mbt| \le 1$,
$$
\Phi_1 (\mbt) = \int_{B_r \setminus U_{\epsilon,r}} |\phi(\mbt + \mbu)| \d \mbu + \int_{B^c_r} |\phi(\mbt + \mbu)| \d \mbu \le \epsilon \| \phi \|_\infty + \int_{B^c_{r-1}} |\phi(\mbu)| \d \mbu,
$$
where we have used $\Leb_\nu (B_r \setminus U_{\epsilon,r}) < \epsilon$. This yields
\begin{equation*}
\delta_{r,\epsilon} :=\sup_{|\mbt| \le 1, \, p=1,2} \Phi_p (\mbt) \to 0 \qquad \text{as } \epsilon \to 0, \, r \to \infty.
\end{equation*}
Next, using
$|\Psi(z)| \le (2|z|)  \wedge (|z|^2/2)$, $z \in \R$, and the Minkowski inequality, we get
\begin{align*}
\int_{U^c_{\epsilon,r}} |\Psi(
\theta \xi_\lambda (\mbu))| \d \mbu
&\le C\int_{U^c_{\epsilon,r}} |\xi_\lambda (\mbu)| \wedge |\xi_\lambda (\mbu)|^2 \d \mbu \\
&\le C\Big(\int_{U^c_{\epsilon,r}} |\xi_\lambda (\mbu)| \d \mbu \Big) \wedge \Big( \int_{U^c_{\epsilon,r}}|\xi_\lambda (\mbu)|^2 \d \mbu \Big) \le C (\zeta_{\lambda,1} \wedge (\zeta_{\lambda,2})^2)
\end{align*}
with
\[
\zeta_{\lambda,p} :=
\lambda^{-\nu/\alpha} \int_{\R^\nu}  \Phi_p (\mbt/\lambda) \I (\mbt  \in \Xi) \d \mbt.
\]
Further, for $p=1$ or $2$,
\begin{eqnarray*}
\zeta_{\lambda,p}
&=& \lambda^{-\nu/\alpha} \Big(\int_{\bar B_\lambda} \Phi_p (\mbt/\lambda) \I (\mbt  \in \Xi)  \d \mbt
+ \int_{\bar B_\lambda^c} \Phi_p (\mbt/\lambda) \I (\mbt  \in \Xi) \d \mbt \Big)\nn \\
&\le&  \lambda^{-\nu/\alpha} (\delta_{r,\epsilon} \Leb_\nu (\Xi) + C \Leb_\nu (\Xi \cap \bar B_\lambda^c) ),
\end{eqnarray*}
where $\bar B_\lambda := \{\mbt \in \R^\nu : |\mbt| \le \lambda \}$ denotes a closed ball of radius $\lambda$ and $\bar B^c_\lambda := \R^\nu \setminus \bar B_\lambda$ denotes its complement.
Hence, denoting
\begin{equation} \label{chicla}
\zeta_\lambda:= \lambda^{-\nu/\alpha} (\delta_{r,\epsilon} \Leb_\nu (\Xi) + C \Leb_\nu (\Xi \cap \bar B_\lambda^c) ),
\end{equation}
the two preceding inequalities give
\begin{equation*}
|j_\lambda(\theta; \phi, U_{\epsilon, r}^c)|\le \lambda^\nu \E \int_{U^c_{\epsilon,r}} |\Psi(\theta \xi_\lambda (\mbu))| \d \mbu
\le C \lambda^\nu \E( \zeta_\lambda \wedge (\zeta_\lambda )^2 ).
\end{equation*}
Now, we shall bound $\E( \zeta_\lambda \wedge (\zeta_\lambda )^2 )$. We can write
\[
\E(\zeta_\lambda\wedge(\zeta_\lambda )^2)
= \E[(\zeta_\lambda)^2 \I (\zeta_\lambda \le 1)]
+ \E [\zeta_\lambda  \I (\zeta_\lambda > 1)]
=\int_0^1  x \P(\zeta_\lambda > x) \d x
+ \int_{1}^\infty  \P(\zeta_\lambda > x) \d x,
\]
where the last equality follows using integration by parts.
Here,
$$
\P ( \zeta_\lambda > x) \le \P ( \zeta^{(1)}_\lambda > x/2) + \P ( \zeta^{(2)}_\lambda > x/2),
$$
where $\zeta^{(1)}_\lambda :=  \lambda^{-\nu/\alpha} \delta_{r,\epsilon} \Leb_\nu (\Xi)$,  $\zeta^{(2)}_\lambda :=
\lambda^{-\nu/\alpha} \Leb_\nu (\Xi \cap \bar B_\lambda^c)$.
By condition \eqref{1volf}, we get
$\P ( \zeta^{(1)}_\lambda > x)
\le C\lambda^{-\nu} \delta_{r,\epsilon}^\alpha x^{-\alpha} $ and therefore
$$
\int_0^1 x \P(\zeta^{(1)}_\lambda > x) \d x
+ \int_{1}^\infty  \P(\zeta^{(1)}_\lambda > x) \d x \le C \lambda^{-\nu} \delta_{r,\epsilon}^{\alpha} \big(\int_0^1 x^{1-\alpha} \d x  +
\int_1^\infty x^{-\alpha} \d x\big)
\le C \lambda^{-\nu} \delta_{r,\epsilon}^{\alpha}.
$$
Similarly, using condition \eqref{2volf},
\begin{eqnarray*}
\int_0^1 x \P(\zeta^{(2)}_\lambda > x) \d x
+ \int_{1}^\infty  \P(\zeta^{(2)}_\lambda > x) \d x
&\le&
\int_0^\infty \P(\zeta^{(2)}_\lambda > x)  \d x =
\E \zeta^{(2)}_\lambda \\
&=&\lambda^{-\nu/\alpha}
o(\lambda^{(\nu/\alpha)(1-\alpha)}) =
o(\lambda^{-\nu}).
\end{eqnarray*}
Therefore, $\E( \zeta_\lambda \wedge (\zeta_\lambda )^2 )\le C\lambda^{-\nu}\delta_{r,\epsilon}^\alpha+o(\lambda^{-\nu})$. It follows that
$|j_\lambda(\theta; \phi, U_{\epsilon, r}^c)|
\le C \delta_{r,\epsilon}^\alpha  + o(1)$, proving
\eqref{jnew}.

Now we consider \eqref{jnew1}. Let $r>0$, $\epsilon>0$ be fixed. Recall that, by Lusin's theorem, there exists a continuous application $\phi_{\epsilon,r} : \R^\nu \to \R$ with compact support in $B_r$ such that $\phi_{\epsilon,r} = \phi$ on $U_{\epsilon,r}$ and moreover, $\| \phi_{\epsilon,r} \|_\infty \le \| \phi \|_\infty$. Since $\Leb_\nu (B_r \setminus U_{\epsilon,r}) < \epsilon$, it suffices to prove \eqref{jnew1} for $\phi_{\epsilon,r}$ on $B_r$ in place of $\phi$ on $U_{\epsilon,r}$. More specifically, it suffices to prove the relation \eqref{jnew1} for integrals over sets $B^+_r := \{\mbu \in B_r: \phi_{\epsilon,r} (\mbu) \ge 0 \}$, $B^-_r := \{\mbu \in B_r: \phi_{\epsilon,r}(\mbu) \le 0 \}$.
Assume without loss of generality that $\phi = \phi_{\epsilon,r} \ge 0$ so that $B^+_r  = B_r$ and $j_\lambda (\theta; \phi, B_r)  = j_\lambda (\theta; \phi)$, \  $j (\theta; \phi, B_r)  = j(\theta; \phi)$.
Then integrating by parts as in \eqref{lalaj},
$$
j_\lambda(\theta; \phi)
= \i \theta \int_{B_r} \Big\{ \int_{\R_+} (\e^{\i  \theta x} -1) \lambda^\nu \P(\xi_\lambda (\mbu) > x) \d x \Big\} \d \mbu.
$$
Hence, $\lim_{\lambda \to \infty} j _\lambda (\theta; \phi)  = j(\theta; \phi)$ follows provided the following two relations hold: for all $x>0$, $\mbu \in B_r$,
\begin{align}\label{limla}
\lim_{\lambda \to \infty} \lambda^\nu x^\alpha \P(\xi_\lambda (\mbu) > x)
&= \ c_\Xi \phi(\mbu)^\alpha, \\
 \text{and }\
 \lambda^\nu x^\alpha \P(\xi_\lambda (\mbu) > x)& < \ C,\nn
\end{align}
with $C$ independent of $\lambda>0$, $x$, and $\mbu$. The second relation in \eqref{limla} follows from \eqref{1volf} since $\lambda^{-\nu/\alpha} \| \phi \|_\infty \Leb_\nu (\Xi) \ge \xi_\lambda (\mbu)$ for all $\mbu \in B_r$. Consider the first one. Note that condition  \eqref{1volf} implies that
$\tilde \xi_\lambda (\mbu) := \lambda^{-\nu/\alpha} \phi(\mbu) \Leb_\nu (\Xi)$ satisfies
$\lim_{\lambda \to \infty} \lambda^\nu x^\alpha \P(\tilde \xi_\lambda (\mbu) > x) = c_\Xi \phi(\mbu)^\alpha $ for all $x >0$, $\mbu \in B_r$. Denote $\eta_\lambda (\mbu) := \xi_\lambda (\mbu) - \tilde \xi_\lambda (\mbu)$. Then for any $\gamma>0$,
\begin{align}\label{xieta}
\P(\xi_\lambda (\mbu) > x)
&\le \P(\tilde \xi_\lambda (\mbu) > (1-\gamma) x)
+ \P(|\eta_\lambda(\mbu)|  >\gamma x), \\
\P(\xi_\lambda (\mbu) > x)
&\ge \P(\tilde \xi_\lambda (\mbu) > (1+\gamma) x) - \P(|\eta_\lambda(\mbu)|  >\gamma x), \nn
\end{align}
It remains to prove that, for any $\gamma>0$,
\begin{eqnarray}\label{limetala}
\lim_{\lambda \to \infty} \lambda^\nu x^\alpha \P(|\eta_\lambda (\mbu)| > \gamma x) = 0.
\end{eqnarray}
Let $\gamma>0$. Since $\phi$ is uniformly
continuous, for any $\delta >0$ there is a $\tau >0$ such  that
$\sup_{|\mbt| \le \tau} |\phi ( \mbt +  \mbu ) - \phi(\mbu)| < \gamma \delta$
uniformly in  $\mbu \in B_r$.
Therefore,
\begin{align*}
|\eta_\lambda(\mbu)| &\le \lambda^{-\nu/\alpha} \int_{\R^\nu} |\phi( \mbt / \lambda + \mbu) -  \phi(\mbu)|  \I (\mbt  \in \Xi)\d \mbt\\
&\le \lambda^{-\nu/\alpha} ( \gamma \delta \Leb_\nu (\Xi) + C \Leb_\nu (\Xi \cap \bar B^c_{\tau \lambda})).
\end{align*}
Thus, the proof of \eqref{limetala} is completely analogous to that
of the estimation of
$\zeta_\lambda$ in \eqref{chicla} and we omit the details.
This also completes the proof of \eqref{jnew1} and thus of \eqref{limvolX}.

\smallskip

\noi {\it Step 2: proof of \eqref{limvolX2}.}
Let $k\ge 1$ and set $G_k(x)=\I(x\ge k)$. We first prove that
\begin{equation}\label{cvX_k}
\lambda^{-\nu/\alpha}(\widehat X_k(\mbt) - \E \widehat X_k(\mbt)) \limd e^{-\mu}\frac{\mu^{k-1}}{(k-1)!}L_\alpha(\phi).
\end{equation}
As noted in introduction, we use the approach in
\cite{sur2024} based on Charlier expansions. Let $N$ be a Poisson
random variable with mean $\E N = \mu$.
Any $G :\N\to\R$ with $\E G(N)^2 < \infty$
can be uniquely expanded in Charlier polynomials $P_j(x;\mu)$, $j =
0,1, \dots$, as
\begin{equation} \label{char0}
G(x) = \sum_{j=0}^\infty \frac{c_{G,\mu}(j)}{j!}P_j(x;\mu)
\end{equation}
with generating function
$$
\sum_{j=0}^\infty \frac{u^j}{j!} P_j(x; \mu)
= (1+ u)^x \e^{-u \mu}, \; u \in \C.
$$
The Charlier coefficients in \eqref{char0} are given by
$ c_{G,\mu}(j) = \mu^{-j} \E[ G(N) P_j(N;\mu)] $ and satisfy
$\sum_{j=0}^\infty  \frac{c^2_{G,\mu}(j)}{j!} \mu^j =
\E G(N)^2 < \infty $. Particularly, for $G_k$ we get
\begin{equation}\label{char}
\widehat X_k(\mbt) - \E \widehat X_k(\mbt) =  G_k(X(t)) - c_{G_k,\mu}(0) = \sum_{j=1}^\infty \frac{c_{G_k, \mu}(j)}{j!}
P_j(X(\mbt); \mu),
\end{equation}
where  $P_1(x;\mu) = x - \mu$ and
\begin{equation*}
c_{G_k,\mu}(1) = \mu^{-1} \E[G_k(N) (N - \mu)]
= \e^{-\mu} \frac{\mu^{k-1}}{(k-1)!} > 0, \quad k =1,2, \ldots.
\end{equation*}
Let
$$
{\mathcal Z}_k(\mbt):=  \widehat X_k(\mbt) - \E \widehat X_k(\mbt) -   c_{G_k,\mu}(1) (X(\mbt) - \E X(\mbt)).
$$
Thus,
$$\widehat X_{\lambda,k}(\phi) -\E \widehat X_{\lambda,k}(\phi)
= c_{G_k,\mu}(1)\left[X_{\lambda}(\phi) -\E X_{\lambda}(\phi)\right] + {\mathcal Z}_{\lambda,k}(\phi),$$
where ${\mathcal Z}_{\lambda,k}(\phi) :=
 \int_{\R^\nu} \phi(\mbt/\lambda) {\mathcal Z}_k(\mbt) \d \mbt $.  Therefore,
\eqref{limvolX2} follows from \eqref{limvolX} provided we can show that
${\mathcal Z}_{\lambda,k}(\phi)$ is negligible, or
\begin{eqnarray}\label{relZ2}
{\rm Var}({\mathcal Z}_{\lambda,k}(\phi))
&=&\int_{\R^{2\nu}} \phi(\mbt_1/\lambda) \phi(\mbt_2/\lambda)
{\rm Cov}({\mathcal Z}_k(\mbt_1),  {\mathcal Z}_k(\mbt_2)) \d \mbt_1 \d \mbt_2 \ = \
o(\lambda^{2\nu/\alpha})
\end{eqnarray}
holds.
To evaluate the last double integral,
we use the following correlation inequality that is a particular case of the last inequality in \cite[Corollary 1]{sur2024}.
\begin{lemma}\label{lem:corr_ineg}
Set $N_1=M_1+M_3$ and $N_2=M_2+M_3$ for $M_1,\,M_2,\,M_3$ independent Poisson random variables with respective means $\mu_1-\mu_3$, $\mu_2-\mu_3$, and $\mu_3$ with $0\le \mu_3<\mu_1\wedge\mu_2$.
Let $G:\N\to\R$ be a function and for $i=1,2$,
\[
k_i^*:=\min\{j\ge1\mid c_{G,\mu_i}(j):=
\mu_i^{-1}\E[G(N_i)P_j(N_i;\mu_i)]\ne0\}.
\]
Then
\begin{equation}\label{char1}
\Cov(G(N_1),G(N_2)) \ =  \
\sum_{j=1}^\infty \frac{c_{G,\mu_1}(j)  c_{G, \mu_2}(j)}{j!} \mu_3^j
\end{equation}
and
\begin{equation}\label{char11}
 |\Cov(G(N_1),G(N_2))|\le \left(\frac{\mu_3}{\sqrt{\mu_1\mu_2}}\right)^{k_1^*\wedge k_2^*} \sqrt{{\rm Var}(G(N_1)){\rm Var}(G(N_2))}.
\end{equation}
\end{lemma}
Note that $\mu_3=\Cov(N_1,N_2)$ and $\Cov(G(N_1),G(N_2)) \ge 0$
for $\mu_1 = \mu_2.$

Applying Lemma 1 to $N_1=X(\boldsymbol{0})$, $N_2=X(\mbt), \, \mu_1 = \mu_2
= {\rm Var}(X(\boldsymbol{0})) = \mu, \, \mu_3= {\rm Cov}(X(\mbt),X(\boldsymbol{0})) =
r_X(\mbt)$, and $G:x\mapsto G_k(x)-c_{G_k,\mu}(1)x, \, k_1^*=k_2^*=2$,
we get 
\begin{eqnarray}\label{relZ1}
|{\rm Cov}({\mathcal Z}_k(\mbt),  {\mathcal Z}_k(\boldsymbol{0}))|
&\le& \mu^{-2} r^2_X(\mbt) {\rm Var} ({\mathcal Z}_k(\boldsymbol{0})) \le
 \mu^{-2} r^2_X(\mbt) {\rm Var} (G_k(\boldsymbol{0}))
\end{eqnarray}
since ${\rm Var}(G(X(\boldsymbol{0})) \le {\rm Var}(G_k(X(\boldsymbol{0}))$.
Observe that
\begin{eqnarray}\label{rXtail}
r_X(\mbt)&=&o  \left(\frac{1}{|\mbt|^{\nu (\alpha-1)/\alpha}}\right), \quad |\mbt| \to \infty.
\end{eqnarray}
Indeed, 
writing $ B_r^c := \R^\nu \setminus B_r$, we have
$$r_X(\mbt)=\E \Leb_\nu (\Xi \cap (\Xi - \mbt)) \
\le\ 2\E \Leb_\nu (\Xi \cap B^c_{|\mbt|/2})
$$
and thus \eqref{rXtail} follows from \eqref{2volf}.
Therefore,
$r_X(\mbt)\le f(|\mbt|)|\mbt|^{\nu (1-\alpha)/\alpha}$ for some bounded continuous function $f$ on $\R_+$ satisfying $\lim_{s\to\infty}f(s)=0$.
By \eqref{relZ1}  and a change of variables,
\begin{align*}
 \lambda^{-2\nu/\alpha}{\rm Var}({\mathcal Z}_{\lambda,k}(\phi))
&\le \mu^{-2}{\rm Var}(G_k(X(\boldsymbol{0})))\int_{\R^{2\nu}} |\phi(\mbt_1) \phi(\mbt_2)| \,|\mbt_1-\mbt_2|^{2\nu(1-\alpha)/\alpha} f^2(\lambda|\mbt_1-\mbt_2|)\d \mbt_1 \d \mbt_2.
\end{align*}
Hence, \eqref{relZ2} follows if we show that the last integral vanishes
with $\lambda \to \infty$. Towards this aim, note that
for all $\mbt_1\ne \mbt_2$, $\phi(\mbt_1) \phi(\mbt_2)|\mbt_1-\mbt_2|^{2\nu(1-\alpha)/\alpha} f^2(\lambda|\mbt_1-\mbt_2|)$ converges to $0$ as $\lambda\to\infty$ and it is dominated by
$\|f\|_\infty^2|\phi(\mbt_1) \phi(\mbt_2)||\mbt_1-\mbt_2|^{2\nu(1-\alpha)/\alpha}$ which is integrable over $\R^{2\nu}$ with respect to $\d\mbt_1 \d\mbt_2$.
Indeed, using that $\phi\in L^1\cap L^\infty$, we obtain
\begin{align*}
 \int_{\R^\nu} \int_{\R^\nu}|\phi(\mbt_1) \phi(\mbt_2)||\mbt_1-\mbt_2|^{2\nu(1-\alpha)/\alpha}\d\mbt_1\d\mbt_2
&\le\int_{\R^\nu} \int_{\R^\nu}|\phi(\mbt_1+\mbt_2) \phi(\mbt_2)|(|\mbt_1|^{2\nu(1-\alpha)/\alpha}\vee 1)\d\mbt_1\d\mbt_2\\
&\le \|\phi\|_\infty\|\phi\|_1\int_{B(0,1)}|\mbt_1|^{2\nu(1-\alpha)/\alpha}\d\mbt_1+\|\phi\|_1^2,
\end{align*}
which is finite since $2\nu(1-\alpha)/\alpha\in (-\nu,0)$.
Hence, \eqref{relZ2} follows by the dominated convergence theorem, thereby
proving \eqref{cvX_k}.

To prove \eqref{limvolX2} we apply the Cram\'er-Wold device. We can use the
previous argument, replacing the function $G_k$ by a linear combination
$G(x)=\sum_{k=1}^na_kG_k(x)$ for some positive integer $n$ and $a_1,\ldots,a_n\in\R$ and noting that
\[
 c_{G,\mu}(1):=\mu^{-1}E[G(N)P_1(N;\mu)]=e^{-\mu}\sum_{k=1}^na_k\frac{\mu^{k-1}}{(k-1)!}.
\]
holds by linearity. We omit the details of this generalization since they are completely analogous to the case $G = G_k$ above.
This completes the proof of Theorem~\ref{thm1}.
\hfill $\Box$

\subsection{Short-range dependence case}

Next, we turn to the central limit theorem (asymptotic Gaussianity of
$X_\lambda (\phi)$) under the short-range dependence condition
\begin{equation} \label{srd}
\E \Leb_\nu (\Xi)^2 < \infty.
\end{equation}

\begin{theorem} \label{thm2} Let $X$ be the random grain model in \eqref{rG1} with generic grain
satisfying $\E \Leb_\nu (\Xi)^2 < \infty$.
Then for any $\phi \in \Phi,$ and  $X_\lambda (\phi)$, $\widehat X_{\lambda, k} (\phi)$, $k\ge 1$, given in \eqref{def:Xlambda}, one has
\begin{eqnarray}\label{limvolX2clt}
\lambda^{-\nu/2}(X_{\lambda}(\phi) -
\E X_{\lambda}(\phi)) &\limd&
W(\phi),  \quad  \mbox{ as }  \lambda \to \infty,
\end{eqnarray}
where $W(\phi)$ is a centered Gaussian variable with variance
$${\rm Var}\left(W(\phi)\right)=\|\phi\|_2^2 \int_{\R^\nu}r_X(\mbt)\d\mbt =\|\phi\|_2^2 \, \E \Leb_\nu (\Xi)^2\in (0,+\infty).$$
Moreover,
\begin{eqnarray}\label{limvolX2fddclt}
\left\{\lambda^{-\nu/2}(\widehat X_{\lambda,k}(\phi) -
\E \widehat X_{\lambda,k}(\phi)); k\ge 1\right\} &\overset{fdd}{\longrightarrow}&
\left\{W_k(\phi), k\ge 1\right\},  \quad  \mbox{ as }  \lambda \to \infty,
\end{eqnarray}
where $\{W_k(\phi), k\ge 1\}$ is a sequence of centered Gaussian variables with  covariances
$$\Cov\left(W_k(\phi),W_l(\phi)\right)=\|\phi\|_2^2\int_{\R^\nu}\Cov\left(\widehat X_{k}(\mbt), \widehat X_{l}(\boldsymbol{0})\right)\d\mbt .$$
\end{theorem}
\begin{remark}\label{rem2}
Using the Charlier decomposition in \eqref{char}, we have
\[
 \Cov\left(\widehat X_{k}(\mbt), \widehat X_{l}(\boldsymbol{0})\right)=\sum_{j=1}^{\infty} \frac{c_{G_k,\mu}(j)c_{G_l,\mu}(j)}{j!} r_X(\mbt)^j.
\]
It follows that
\[
 \Cov\left(W_k(\phi),W_l(\phi)\right)=c_{G_k,\mu}(1)c_{G_l,\mu}(1)\Var\left(W(\phi)\right)+\|\phi\|_2^2\sum_{j=2}^\infty \frac{c_{G_k,\mu}(j)c_{G_l,\mu}(j)}{j!} \int_{\R^\nu}r_X(\mbt)^j \d\mbt.
\]
This is in contrast with Remark~\ref{rem1} since all Charlier coefficients of
\eqref{char} contribute to the limit distribution.
\end{remark}

Note that, as a corollary we recover Theorem~3.5 of \cite{Hall88}, that we state here for completeness.
\begin{corollary}\label{Cor:thm2} Under the assumptions of Theorem \ref{thm2},
 for any bounded Borel set
 $A \subset \R^\nu$ with $\Leb_\nu (A) >0$,
\begin{equation*}
\lambda^{\nu/2}(\widehat p_{\lambda, A} -
p) \limd
\sigma W(A)/\Leb_\nu(A),  \quad  \mbox{ as } \lambda \to \infty
\end{equation*}
where $W(A) \sim N(0, \Leb_\nu(A))$ and
$$
\sigma^2 :=   \int_{\R^\nu} {\rm Cov}(\widehat X(\boldsymbol{0}), \widehat X(\mbt)) \d \mbt
= \e^{-2\mu} \int_{\R^\nu}
\big(\e^{\Leb_\nu (\Xi \cap (\Xi - \mbt))} -1 \big) \d \mbt.
$$
\end{corollary}

\noi {\bf Proof of Theorem~\ref{thm2}.}
Let $m \ge 1$. To show \eqref{limvolX2clt},  we can use an approximation
by a $m$-dependent random grain random field and the central limit theorem for such random fields.
Consider
\begin{eqnarray*}
X^{(m)}(\mbt) = \sum_{j=1}^\infty \I (\mbt \in (\mbu_j + \Xi_j \cap B_{m/2})), \quad \mbt \in  \R^\nu
\end{eqnarray*}
the random grain model with generic grain $\Xi \cap B_{m/2} \subset \{\mbt\in\R^\nu;|\mbt| \le m/2\}$ belonging to the ball $B_{m/2}$ of radius $m/2 $.
Thus, $X^{(m)} (\mbt_1) $ and $X^{(m)}(\mbt_2)$
are independent when $|\mbt_1 - \mbt_2| > m $. This fact
follows from the independence property of Poisson stochastic integrals with disjoint  supports, as
$$\I (\mbt_i - \mbu \in \Xi \cap B_{m/2}, i=1,2)
\le \I (|\mbt_i - \mbu| \le m/2, i=1,2)
\le \I ( |\mbt_1 - \mbt_2| \le m) = 0 $$  for any $\mbu \in \R^\nu$ by triangle inequality. Note that
\begin{eqnarray}
r_{X^{(m)}}(\mbt)
&:=&{\rm Cov}(X^{(m)}(\boldsymbol{0}),
X^{(m)} (\mbt)) = \E \Leb_\nu ((\Xi \cap B_{m/2}) \cap [(\Xi \cap B_{m/2}) - \mbt]) \nn \\
&\le&\E \Leb_\nu (\Xi \cap (\Xi - \mbt)) = r_{X}(\mbt)
\label{rrm}
\end{eqnarray}
and $ r_{X^{(m)}}(\mbt) \nearrow r_{X}(\mbt) $ as $m \to \infty$
at each $\mbt$.
Similarly observe that
\[
 |\Cov(X(\mbt)-X^{(m)}(\mbt),X(\boldsymbol{0})-X^{(m)}(\boldsymbol{0}))|
 =\E\Leb_\nu((\Xi\cap B_{m/2}^c)\cap(\Xi\cap B_{m/2}^c-\mbt)).
\]
This last expression converges to $0$ as $m\to\infty$ and it is uniformly bounded by $r_X(\mbt)$.
Defining as before, $X^{(m)}_\lambda(\phi)=\int_{\R^\nu}\phi(\mbt/\lambda)X^{(m)}(\mbt)\d\mbt$, by the dominated convergence theorem, we get that
\begin{eqnarray*}
&&\hspace{-30pt}\lambda^{-\nu} {\rm Var}( X_{\lambda}(\phi) -  X_{\lambda}^{(m)}(\phi)) \\
&=&\lambda^{-\nu} \int_{\R^{2\nu}} \phi(\mbt_1/\lambda) \phi(\mbt_2/\lambda) \,
{\rm Cov}( X(\mbt_1) -  X^{(m)}(\mbt_1), X(\mbt_2) -  X^{(m)}(\mbt_2)) \d \mbt_1 \d \mbt_2 \\
&\le&C \lambda^{-\nu} \int_{\R^{\nu}} |\phi(\mbt_1/\lambda)| \d \mbt_1 \times
\int_{\R^{\nu}}|{\rm Cov}
( X(\boldsymbol{0}) -  X^{(m)}(\boldsymbol{0}), X(\mbt_2) -  X^{(m)}(\mbt_2))|  \d \mbt_2
\ \to  \ 0,
\end{eqnarray*}
as $m \to \infty$, uniformly in $\lambda>0$.
The convergence in
\eqref{limvolX2clt} then follows from the central limit theorem for $m$-dependent random fields
\[
\lambda^{-\nu/2} ( X^{(m)}_{\lambda}(\phi) -
\E  X^{(m)}_{\lambda}(\phi))
\limd  W^{(m)}(\phi),
\]
where $W^{(m)}(\phi)$ is a centered Gaussian random variable with
$$
\Var(W^{(m)}(\phi))= \|\phi\|_2^2
\int_{\R^\nu} {\rm Cov}(X^{(m)}(\boldsymbol{0}), X^{(m)}(\mbt)) \d \mbt \to \|\phi\|_2^2\int_{\R^\nu}r_X(\mbt)\d\mbt  < \infty
$$
as $m \to \infty$.

\medskip

We prove \eqref{limvolX2fddclt} in a similar way, using Charlier expansion. We first define, for some $k\ge 1$, the approximating random field $\widehat X^{(m)}_k$ as $\widehat X^{(m)}_k(\mbt) := G_k(X^{(m)}(\mbt))$ with $G_k(x)=\I(x \ge k)$.
Applying Lemma~\ref{lem:corr_ineg} with $N_1=X^{(m)}(\boldsymbol{0})$, $N_2=X^{(m)}(\mbt)$, and $G=G_k$, and using \eqref{rrm},
we have that
\begin{eqnarray*}
|{\rm Cov}(\widehat X^{(m)}_k(\boldsymbol{0}),
\widehat X^{(m)}_k (\mbt))|
&\le&\frac{r_{X^{(m)}}(\mbt)}{r_{X^{(m)}}(\boldsymbol{0}) }\,
{\rm Var}(\widehat X^{(m)}_k(\boldsymbol{0})) \ \le \ C r_X(\mbt).
\end{eqnarray*}
In a similar way,
\begin{eqnarray*}
|{\rm Cov}(\widehat X_k(\boldsymbol{0})- \widehat X^{(m)}_k(\boldsymbol{0}),
\widehat X_k(\mbt)- \widehat X^{(m)}_k (\mbt))|
&\le&C r_X(\mbt)
\end{eqnarray*}
is bounded by an integrable function uniformly in $m \ge 1
$ and vanishes with $m \to \infty$ at each point
$\mbt \in \R^\nu$.
Arguing as before we get that $\lambda^{-\nu}\Var( \widehat X_{\lambda,k}(\phi) -  \widehat X_{\lambda,k}^{(m)}(\phi))\to 0$
as $m\to\infty$ uniformly in $\lambda$. We can thus deduce the CLT for $\widehat X_{\lambda,k}(\phi)$ from the CLT for the $m$-dependent random field $(\widehat X^{(m)}_{\lambda,k}(\mbt))$.

The finite dimensional convergence can be obtained similarly by applying the Cram\'er-Wold device. Again, we omit the details.
\hfill $\Box$

\medskip

\begin{remark}\label{rem:srd}
In \cite{kulik2021long},  the authors introduced another notion of short-range dependence
which applies to random fields with infinite variance and involves
the covariance function of excursion sets.
More precisely, they call  a stationary field $Y$  short-range dependent if
\begin{equation}\label{newSRD}
\int_{\R^\nu}\int_{\R}\int_{\R}{\rm Cov}(\I(Y (\boldsymbol{0})>u), \I(Y (\mbt)>v)) \d F(u) \d F(v) \d\mbt \, < \, \infty
\end{equation}
holds for all distribution functions $F$ on $\R$. A sufficient condition for infinitely divisible fields to satisfy \eqref{newSRD} is also given in \cite{makogin2025short}.
Observe that for our Poisson random grain model in \eqref{rG1},
 \eqref{newSRD} is equivalent to \eqref{srd},
 that is, to the short-range dependence as defined in introduction and used in Theorem 2.
Let us prove this equivalence.
%
First, using that
$
\int_{\R}
\I(X(\mbt) > u) \d F(u) = G(X(\mbt))$ for $G(x):=F(x^-)$ the left-continuous version of $F$,
we see that \eqref{newSRD} for $Y =X$ writes as
\begin{equation*}
\int_{\R^\nu} {\rm Cov}(G(X (\boldsymbol{0})), G(X (\mbt)))  \d\mbt \, < \, \infty.
\end{equation*}
According to Charlier expansion in \eqref{char1}
\begin{equation*}
\Cov(G(X(\boldsymbol{0})),G(X(\mbt))) \ =  \
\sum_{j=1}^\infty \frac{c^2_{G,\mu}(j)}{j!} r^j_X(\mbt).
\end{equation*}
Since $r_X(\mbt) \ge 0$, we get
\begin{equation*}
\int_{\R^\nu} {\rm Cov}(G(X (\boldsymbol{0})), G(X (\mbt)))  \d\mbt \, \ge \, c^2_{G,\mu}(1)  \int_{\R^\nu} r_X(\mbt) \d\mbt.
\end{equation*}
Since $c_{G,\mu}(1) = \E [G(X(\mbt)+1) - G(X(\mbt))] >0$ as soon as there exists $x>0$ such that $G(x)<1$,
we see that  \eqref{newSRD} implies
$\int_{\R^\nu} r_X(\mbt) \d\mbt < \infty$. Conversely,
from \eqref{char11},  $|\Cov(G(X(\boldsymbol{0})),G(X(\mbt)))|
\le r_X(\mbt) {\rm Var}(G(X(\boldsymbol{0})))/\mu $, implying
\begin{equation*}
\int_{\R^\nu} \Cov(G(X(\boldsymbol{0})),G(X(\mbt))) \d \mbt
 \le   ({\rm Var}(G(X(\boldsymbol{0})))/\mu)
\int_{\R^\nu } r_X(\mbt) \d \mbt
\end{equation*}
and completing the proof of the converse implication.
\end{remark}

\subsection{Discussion}\label{sec:disc}

As noted in several works, limit theorems for random grain models are related with aggregated network traffic, particularly,
with scaling limits of aggregated workload in standard M/G/$\infty$
queuing systems and related Poisson shot-noise processes \cite{bier2010, kaj2007, kajt2008, leip2023, miko2002,
pils2016, sur2024}. The (aggregated) infinite source Poisson model can  be defined as a Poisson stochastic integral
$$
Y_M(t) = \int_{\R \times \R_+^2} a \I (u < t \le u +r) {\cal N}_M (\d u, \d a, \d r), \quad t \in \R,
$$
with respect to Poisson random measure
${\cal N}_M (\d u, \d a, \d r)$ so that a Poisson event in $(u, a, r)$ represents a session arriving at time $u$ of duration $r$ and load $a$. The mean
$E {\cal N}_M (\d u, \d a, \d r) = M \d u  \P( (A,R)
\in (\d a, \d r) )$ is proportional to the
parameter $M$  called the {\it connection rate}. The problem studied
in the above mentioned and other related works concerns the limit distribution of the {\it aggregated workload}
\begin{equation*} 
\mathcal A(\lambda y, M) := \int_0^{\lambda y} Y_M(t) \d t, \quad y >0,
\end{equation*}
as $\lambda \to \infty$ and $M$ increases with $\lambda$ at a certain
rate. From elementary properties of Poisson integral, $\mathcal A(\lambda y, M)$
can be approximated by the integral
\begin{equation}\label{XM}
X_\lambda(\phi_{M,\lambda}) = \int_{\R^2} \phi_{M,\lambda}(\mbt/\lambda) X(\mbt)\d\mbt
=\int_{ (0,M]\times (0,\lambda y]  } X(t_1, t_2)
\d t_1 \d t_2
\end{equation}
of the random grain model $X(\mbt)$, $\mbt = (t_1,t_2)\in \R^2$ in
\eqref{rG1} with
generic grain
\begin{equation*}
\Xi = [0, A] \times [0,R] \subset \R^2
\end{equation*}
and $\phi_{M,\lambda}(t_1,t_2) = \I((t_1,t_2) \in [0,M/\lambda] \times [0,y])$.
 Limit distributions of
$\mathcal A(\lambda y, M)$  and $X_\lambda(\phi_{M,\lambda})$
as $\lambda \to \infty$, $M = \lambda^\gamma \to \infty$ for some $\gamma >0$
were studied in  \cite{kajt2008,  miko2002, pils2016} and elsewhere under various assumptions on the joint distribution of $(A, R)$.

Note that for
$M = \lambda$,  $\phi_{M,\lambda}(t_1,t_2)=\phi(t_1,t_2) = \I((t_1,t_2) \in [0,1] \times [0,y])$ does not depend on the parameters $M$ and $\lambda$ and the integral in \eqref{XM}
agrees with
$X_\lambda (\phi)$ in Theorem~\ref{thm1}.
These aforementioned works may explain the necessity of
condition \eqref{2volf} for Theorem~\ref{thm1} and
show that the scaling behavior without it
can be non-typical and quite complex, as illustrated in the following examples.


\begin{example} {\rm Let
\begin{eqnarray}\label{Xiex1}
\Xi = [0,1] \times [0, R] \subset \R^2
\end{eqnarray}
where $R>0$ is a random variable with distribution as in Proposition \ref{corRGhomo} (i). 
 According to  \cite{kajt2008, leip2023, miko2002} and other
works, for indicator functions $\phi(\mbt) = \I(\mbt \in [\boldsymbol{0}, \mbx])$, $\mbx \in \R^2_+$,
$\lambda^{-1}X_\lambda (\phi)$ is asymptotically normal, hence does not satisfy
\eqref{limvolX}. It is easy to see that  \eqref{Xiex1} violates condition  \eqref{2volf}: for large $\lambda$ we have
$\E  \Leb_2 \big(\Xi \cap   \{|\mbt| > \lambda\}\big) \sim \E (R - \lambda)_+ = O(\lambda^{-(\alpha -1)}) \gg   O(\lambda^{-2(\alpha -1)/\alpha}) $
since $\alpha < 2 $. }

\end{example}

\begin{example} (Deterministically related transmission rate and duration model, see \cite{pils2016, leip2023} ) {\rm Let
\begin{eqnarray}\label{Xiex2}
\Xi = [0,R^{1-p}] \times [0, R^p] \subset \R^2
\end{eqnarray}
where $p \in (0,1)$ is a (shape) parameter and
 $R>0$ the same as in  \eqref{Xiex2}. By symmetry, it suffices to consider the case $p \ge 1/2 $.
 According to these works, the limit distribution of $X_\lambda (\phi)$ for indicator functions $\phi$
is Gaussian if $\alpha > 2p$ and
$\alpha$-stable if $\alpha < 2p$; for $\alpha = 2p$ this limit is an intermediate one written as a Poisson stochastic integral.
We have $\Leb_2(\Xi) = R$ and
$\Leb_2 \big(\Xi \cap   \{|\mbt| > \lambda\}\big)  < 2 (R^p - \lambda)R^{1-p} $ so that
$\E  \Leb_2 \big(\Xi \cap   \{|\mbt| > \lambda\}\big) = O(\lambda^{-(\alpha -1)/p}) =
o(\lambda^{-(\alpha-1)(2/\alpha)}) $
when $2p < \alpha $ and
 \eqref{2volf} holds. Similarly, one can check
that for  $2p \ge \alpha $   condition  \eqref{2volf} is not satisfied and
Theorem 1 for \eqref{Xiex2} does not apply.
 }

\end{example}

\section{Randomly homothetic Boolean set}\label{sec:homothetic}

In this section we introduce a specific class of random closed sets and give some examples.

\begin{definition} {\rm A random closed set $\Xi \subset \R^\nu$  is called a {\it random homothetic grain} if it can be represented as
\begin{eqnarray}\label{Theta1}
\Xi = R^{1/\nu} \, \Xi^0,
\end{eqnarray}
where $\Xi^0 $ is a random closed set such that  $\Xi^0 \subset \bar{B}_1:=\{\mbt\in\R^\nu;|\mbt| \le 1\} $ a.s.\ and   
$\Leb_\nu (\Xi^0) >0$ a.s., and
$R>0$ is a random variable, independent of $\Xi^0$.
}
\end{definition}

\begin{proposition} \label{corRGhomo} Let $\Xi$ be a random
homothetic grain in \eqref{Theta1} for a positive random variable $R$ and a random closed set $\Xi^0$.

\noi  (i) Assume that there exist $\alpha \in (1,2)$, $c_R >0$ such that
$$\P(R > x) \sim c_R x^{-\alpha},\;  x \to  \infty  .$$  Then
$\Xi$ satisfies
conditions \eqref{1volf} and \eqref{2volf}, with $c_\Xi = c_R \E \Leb_\nu (\Xi^0)^\alpha $. Moreover,
the covariance $r_X(\mbt) = {\rm Cov}(X (\boldsymbol{0}), X(\mbt)) $ in \eqref{covRG} satisfies
\begin{eqnarray} \label{covXRG}
r_X(\mbt)
&=&O(|\mbt|^{-\nu(\alpha -1)}),  \quad |\mbt| \to \infty.
\end{eqnarray}

\noi (ii) Assume that the random variable $R$ has a density $f$ and that  there exist $\alpha \in (1,2)$, $c_f >0$ such that
\begin{equation} \label{falpha}
f(r) \sim c_f r^{-1-\alpha}, \qquad r \to \infty.
\end{equation}
and assume also that the function $\mbt \mapsto \E \Leb_\nu(\Xi^0 \cap (\Xi^0- \mbt)) $ is continuous on $\R^\nu \setminus \{\boldsymbol{0}\}$.
Then
\begin{eqnarray*}
r_X(\mbt)&=&|\mbt|^{-\nu(\alpha -1)} \big(\ell(\frac{\mbt}{|\mbt|}) + o(1)\big),  \quad |\mbt| \to \infty,
\end{eqnarray*}
with
$\ell(\cdot) \in S_{\nu-1}$ on the unit sphere $S_{\nu-1}$ of $\R^\nu $
given by
\begin{eqnarray*}
\ell(\mbz)&:=&c_f \int_{\R_+} \E \Leb_\nu \big(\Xi^0 \cap (\Xi^0 - r^{-1/\nu} \, \mbz)\big) r^{-\alpha} \d r.
\end{eqnarray*}

\smallskip

\noi (iii) Assume that $\E R^2 < \infty $.  Then $\E \Leb_\nu (\Xi)^2 < \infty $.

\end{proposition}

\noi  {\bf Proof.} (i) Condition  \eqref{1volf} follows from $\Leb_\nu(\Xi) = R \Leb_\nu (\Xi^0)$ and Breiman's lemma (see \cite{breiman1965},
\cite[p.231-232]{resnick2007},
\cite[Lemma 1.4.3]{kulikSoulier2020}).
Consider  \eqref{2volf}.
Since $\Xi^0 \subset \bar B_1 $ is bounded, 
$$
\Leb_\nu \big(\Xi \cap   \{|\mbt| > \lambda \}\big)  \le
\Leb_\nu \big( \{\lambda < |\mbt| \le R^{1/\nu} \} \big)
\le C (R - \lambda^\nu) \vee 0
$$
and therefore
\begin{eqnarray*}
\E \Leb_\nu \big(\Xi \cap   \{|\mbt| > \lambda\}\big)
&\le&C \int_{\lambda^\nu}^\infty \P(R> r) \d r   \ = \  O(\lambda^{-\nu(\alpha -1)}) = o(\lambda^{-\nu(\alpha -1)/\alpha})
\end{eqnarray*}
since $\alpha > 1 $.
Part (ii) is similar to \cite[Proposition 1]{sur2024}.
Part (iii) is obvious.  \hfill $\Box$

This notion of random homothetic grain allows to go beyond the classical case of deterministic shapes. The following examples illustrate this fact by giving some possible random shapes.

\begin{example} \label{exhard} (Hard balls grain.) {\rm  Let $\{\mbu^0_j\}$ be a Poisson process on $B_1 \subset \R^\nu$,  with Lebesgue intensity. From each
point $\mbu^0_j$ a hard closed ball
starts growing with unit rate and the growth stops after it hits another ball or the boundary  $\partial B_1 = \{|\mbu| = 1 \}$.
The set $\Xi^0$ is defined as the union of all such balls.
}
\end{example}

\begin{figure}
\centering
        \includegraphics[width=0.35\textwidth]{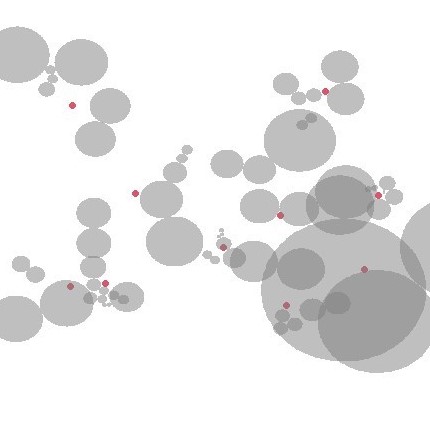} ~~~~~
        \includegraphics[width=0.35\textwidth]{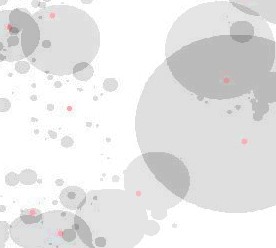}
    \caption{\small Sample of random grain models with homothetic grain sets given by the shapes $\Xi^0$ in Example~\ref{exhard} [left] and Example~\ref{exBool} [right].  The red points correspond to the Poisson  process of germs in \eqref{rG1}.}
    \label{fig:my_label}
\end{figure}

\begin{remark} \rm The hard balls property in Example 1 refers to
individual grain sets, whereas the corresponding Boolean model does not
have this property and allows for overlapping grains as seen in Figure~\ref{fig:my_label}, left. Hard-core random sets with long-range dependence can be obtained from  the Boolean model
by thinning operations, see \cite{matern1960} and more recent work
\cite{KurLes13}.
\end{remark}

\begin{example} \label{exBool} (Cluster Boolean grain.) {\rm Let $\{(\mbu^0_j, y^0_j)\}$ be a Poisson process on $B_1\times (0,1] \subset \R^\nu \times \R_+$ with intensity
$\mu^0(\d \mbu, \d y)$ satisfying  $\int_{B_1 \times (0,1]} y \mu^0 (\d \mbu, \d y) < \infty $, and
$$
\tilde\Xi^0 = \bigcup_{j=1}^\infty (\mbu^0_j + y_j^{1/\nu} B_1).
$$
Then $\tilde\Xi^0 \subset B_2 = \{|\mbu|  < 2\} $ is a.s.\ bounded and is a union of an infinite number of balls unless the intensity measure
$\mu^0(B_1 \times (0,1]) < \infty $  is bounded.
Then, one can choose $\Xi^0:=\frac12\overline{\tilde\Xi^0}\subset\bar B_1$, where $\overline{A}$ denotes the closure of the set $A$.}

\end{example}

\section{Sample volume fraction on hyperplane}\label{sec:hyperplan}

Let $\nu\ge 2$, $\nu_0\in\{1,\ldots,\nu-1\}$,  and define
\begin{eqnarray}\label{hyp0}
H_{\nu_0}&:=&\{\mbt=(t_1,\ldots,t_\nu) \in \R^\nu: t_i = 0, \nu_0 < i \le \nu\}.
\end{eqnarray}
The above hyperspace can be identified with $\R^{\nu_0}$.
We use the notation $\mbt = (\mbt', \mbt'')$, with
$\mbt' := (t_1, \cdots, t_{\nu_0}) \in \R^{\nu_0}$, $\mbt'' =  (t_{\nu_0 +1}, \cdots, t_\nu) \in \R^{\nu- \nu_0} $.

We introduce the integrals of the random grain models observed from the hyperspace $H_{\nu_0}$. For $\phi \in \Phi_0 := L^1(\R^{\nu_0}) \cap L^\infty (\R^{\nu_0})$, set
\[
 X_{0,\lambda}
(\phi) := \int_{\R^{\nu_0}} \phi(\mbt'/\lambda) X(\mbt', 0, \cdots, 0) \d \mbt'
\]
and for $k\ge 1$,
\[
\widehat X_{0,\lambda,k}
(\phi) := \int_{\R^{\nu_0}} \phi(\mbt'/\lambda) \widehat X_k(\mbt', 0, \cdots, 0) \d \mbt',
\]
with  $\widehat X_k$ defined in \eqref{Xk}.

\begin{theorem} \label{thm3} 
Let $X$ be the random grain model in \eqref{rG1} with
randomly homothetic grain set as in \eqref{Theta1} satisfying
the conditions of Proposition \ref{corRGhomo} (i) for some $\alpha\in (1,2)$ and $\nu_0\in \{1,\ldots,\nu-1\}$.

\smallskip

\noi (i) If $1< \alpha < 1+ \frac{\nu_0}{\nu}$, let
\begin{equation*}
\alpha_0 = 1 + \frac{\nu}{\nu_0}(\alpha -1) \in (1,2).
\end{equation*}
Then
\begin{eqnarray}\label{limvolX0}
\lambda^{-\nu_0/\alpha_0}(X_{0,\lambda}(\phi) -
\E X_{0,\lambda}(\phi)) &\limd&
L_{\alpha_0}(\phi),\qquad \mbox{ as } \lambda \to \infty,
\end{eqnarray}
where $L_{\alpha_0}(\phi)$ has an $\alpha_0$-stable distribution with log-characteristic function as in \eqref{limjRG} replacing $\alpha$ by $\alpha_0$, $\R^\nu$  by $\R^{\nu_0}$, and $c_\Xi$ by $c_{0,\Xi}$ defined at \eqref{c0Xi}.
Further,
\begin{eqnarray}\label{limvolXk}
\left\{\lambda^{-\nu_0/\alpha_0}(\widehat X_{0,\lambda,k}(\phi) -
\E \widehat X_{0,\lambda,k}(\phi));\,k\ge 1\right\} &\limfdd&
\left\{e^{-\mu}\frac{\mu^{k-1}}{(k-1)!}L_{\alpha_0}(\phi);\, k\ge 1\right\},
\end{eqnarray}
as $\lambda \to \infty$.

\smallskip

\noi (ii)  If $1+ \frac{\nu_0}{\nu}< \alpha <2$, then
\begin{eqnarray*}
\lambda^{-\nu_0/2}(X_{0,\lambda}(\phi) -
\E X_{0,\lambda}(\phi)) &\limd&
W_{0}(\phi),\qquad  \mbox{ as } \lambda \to \infty,
\end{eqnarray*}
where $W_{0}(\phi)$ is a centered Gaussian random variable with variance
\[
 \Var(W_0(\phi))=\|\phi\|_2^2\int_{\R^{\nu_0}}r_X(\mbt',0,\ldots,0)\d\mbt'.
\]
Further
\begin{eqnarray*}
\left\{\lambda^{-\nu_0/2}(\widehat X_{0,\lambda,k}(\phi) -
\E \widehat X_{0,\lambda,k}(\phi));\,k\ge 1\right\} &\limfdd&
\left\{W_{0,k}(\phi);\, k\ge 1\right\},\qquad  \mbox{ as } \lambda \to \infty,
\end{eqnarray*}
where $\left\{W_{0,k}(\phi);\, k\ge 1\right\}$ is a sequence of centered Gaussian random variables with covariances
\[
 \Cov\left(W_{0,k}(\phi),W_{0,l}(\phi)\right)=\|\phi\|_2^2\int_{\R^{\nu_0}}\Cov\left(\widehat X_{k}(\mbt',0,\ldots,0), \widehat X_{l}(\boldsymbol{0})\right)\d\mbt'.
\]
\end{theorem}
\begin{remark}
The comments of Remark~\ref{rem1} and Remark~\ref{rem2} apply similarly here.
\end{remark}

As a particular case we obtain the following corollary for volume fraction estimator on an hyperplane.

\begin{corollary} \label{cor3} 
Under the assumptions of Theorem \ref{thm3}, and
$\widehat p_{\lambda, A} (H_{\nu_0}) $
in \eqref{volumefgamma} be the volume fraction estimator
on the hyperspace \eqref{hyp0}, where $A \subset H_{\nu_0}$ is an arbitrary bounded
Borel set with $\Leb_{\nu_0} (A) >0$.

\smallskip

\noi (i) If $1< \alpha < 1+ \frac{\nu_0}{\nu}$ and $\alpha_0 = 1 + \frac{\nu}{\nu_0}(\alpha -1)$,
then
\begin{eqnarray*}
\lambda^{\nu_0 -(\nu_0/\alpha_0)}(\widehat p_{\lambda, A}(H_{\nu_0}) -p) &\limd&
\e^{-\mu} L_{\alpha_0}(A)/\Leb_{\nu_0}(A),  \qquad  \mbox{ as } \lambda \to \infty,
\end{eqnarray*}
where $ L_{\alpha_0}(A): = L_{\alpha_0}(\I(A))$ with $L_{\alpha_0}$ given in Theorem~\ref{thm3}.

\smallskip

\noi (ii)  If $1+ \frac{\nu_0}{\nu}<\alpha <2$, then
\begin{eqnarray*}
\lambda^{\nu_0/2}(\widehat p_{\lambda, A}(H_{\nu_0}) -p) &\limd&
W_{0,1}(A)/\Leb_{\nu_0}(A),  \qquad  \mbox{ as }\lambda \to \infty,
\end{eqnarray*}
where $W_{0,1}(A):=W_{0,1}(\I(A))$ for $W_{0,1}$ given in Theorem~\ref{thm3}.
\end{corollary}

\noi {\bf Proof of Theorem~\ref{thm3}.}
(i) We proceed similarly as in the proof of Theorem~\ref{thm1} and first prove \eqref{limvolX0}.
Following \eqref{jlambda},
\begin{eqnarray*}
j_{0,\lambda}(\theta)
&:=&\int_{\R^{\nu_0}\times \R^{\nu-\nu_{0}}} \E \Psi \Big(\frac{\theta}{\lambda^{\nu_0/\alpha_0}}
\int_{\R^{\nu_0}}\phi(\mbt'/\lambda) \I ((\mbt',\boldsymbol{0}) - (\mbs',\mbs'') \in \Xi) \d \mbt' \Big) \d \mbs' \d \mbs'' \\
&=&\lambda^{\nu_0} \int_{\R^{\nu_0} \times \R^{\nu-\nu_{0}}} \E \Psi \Big
(\frac{\theta}{\lambda^{\nu_0/\alpha_0}} \int_{\R^{\nu_0}}
\phi\big(\frac{\mbt'}{\lambda} +  \mbs'\big) \I ((\mbt',\boldsymbol{0}) \in (\boldsymbol{0},\mbs'')  + \Xi) \d \mbt' \Big)
\d \mbs' \d \mbs'', \nn
\end{eqnarray*}
with $\Psi(z)=e^{iz}-1-iz$.
The intuitive argument leading to $j_{0,\lambda} (\theta)  \to j_0 (\theta) $
uses the observation that, as $\lambda\to\infty$,
\begin{eqnarray*}
\int_{\R^{\nu_0}}\phi\big(\frac{\mbt'}{\lambda} +  \mbs'\big) \I ((\mbt',\boldsymbol{0}) \in (\boldsymbol{0},\mbs'') + \Xi) \d \mbt'
&\to&\phi(\mbs')\int_{\R^{\nu_0}} \I ((\mbt',\boldsymbol{0}) \in (\boldsymbol{0},\mbs'')  + \Xi) \d \mbt' \\
&=&\phi(\mbs') \Leb_{\nu_0}(\Xi_{\mbs''}), \nn
\end{eqnarray*}
where $\Xi_{\mbs''} := \Xi \cap  \{\mbt'' = \mbs''  \} $ is the section of $\Xi$ by the hyperplane  $\{\mbt =  (\mbt', \mbt'') \in \R^\nu: \mbt'' = \mbs''  \} $.
As in the proof of Theorem~\ref{thm1} we first consider the limit of
\begin{eqnarray*}
\tilde j_{0,\lambda}(\theta)
&:=&\lambda^{\nu_0}
\int_{\R^{\nu_0} \times  \R^{\nu-\nu_{0}}} \E \Psi \Big(
\frac{\theta}{\lambda^{\nu_0/\alpha_0}} \phi(\mbs') \Leb_{\nu_0}(\Xi_{\mbs''})\Big) \d \mbs' \d \mbs''
\end{eqnarray*}
and then show that the difference $j_{0,\lambda}(\theta) - \tilde j_{0,\lambda}(\theta)$ is negligible.
Since
$\Xi = R^{1/\nu} \, \Xi^0  $ we get  $\Leb_{\nu_0}(\Xi_{\mbs''})
= R^{\nu_0/\nu} g^0(\mbs''/R^{1/\nu}), $ where
\begin{equation}\label{g0}
g^0(\mbs'') := \Leb_{\nu_0}(\Xi^0 \cap \{\mbt''
= \mbs''\}) \ge 0
\end{equation}
is  the $\nu_0$-dimensional Lebesgue measure of the intersection of $\Xi^0$ with the hyperplane
$\{(\mbt',\mbt'') \in \R^\nu: \mbt''= \mbs''\}$. Since $\Xi^0 \subset \{|\mbt| \le 1\} $ is a bounded set,
$g^0(\mbs'')$ in \eqref{g0} is bounded and has a bounded support (vanishes for $|\mbs''| > 1$). Moreover,
$g^0(\mbs'')$ is {independent} of $R$.

Then, similarly as we did in the proof of Theorem \ref{thm1}, using integration by parts we get
\begin{eqnarray*}
 \Psi \Big(
\frac{\theta}{\lambda^{\nu_0/\alpha_0}} \phi(\mbs') \Leb_{\nu_0}(\Xi_{\mbs''})\Big)&=& \i \theta \phi(\mbs')
\int_{\R^+} \Big(\e^{\i \theta\phi(\mbs') x} - 1\Big) \I( \Leb_{\nu_0}(\Xi_{\mbs''}) > x \lambda^{\nu_0/\alpha_0})\d x\\
&=&\i \theta \phi(\mbs')
\int_{\R^+} \Big(\e^{\i \theta\phi(\mbs') x} - 1\Big) \I( R^{\nu_0/\nu}g^0(\mbs''/R^{1/\nu}) > x \lambda^{\nu_0/\alpha_0})\d x
\end{eqnarray*}
Then,
$$\tilde j_{0,\lambda}(\theta)
=\lambda^{\nu_0}\int_{\R^{\nu_0}\times \R^{\nu-\nu_0}}\E\left(\i \theta \phi(\mbs')
\int_{\R^+} \Big(\e^{\i \theta\phi(\mbs') x} - 1\Big) \I( R^{\nu_0/\nu}g^0(\mbs''/R^{1/\nu}) > x \lambda^{\nu_0/\alpha_0})\d x \right)\d \mbs' \d \mbs''$$
and Fubini's Theorem and a change of variables give
\begin{align*}
&\tilde j_{0,\lambda}(\theta)\\
&=\lambda^{\nu_0}\int_{\R^{\nu_0}}\i \theta \phi(\mbs')\E\left( \int_{\R^{\nu-\nu_0}}
\int_{\R^+} \Big(\e^{\i \theta\phi(\mbs') x} - 1\Big) R^{1-\nu_0/\nu} \I( R^{\nu_0/\nu}g^0(\mbs'') > x \lambda^{\nu_0/\alpha_0})\d x \d \mbs'' \right)\d \mbs'\\
&=\lambda^{\nu_0}\int_{\R^{\nu_0}}\i \theta \phi(\mbs') \int_{\R^{\nu-\nu_0}}
\int_{\R^+} \Big(\e^{\i \theta\phi(\mbs') x} - 1\Big) \E\left(R^{1-\nu_0/\nu} \I( Rg^0(\mbs'')^{\nu/\nu_0} > x^{\nu/\nu_0} \lambda^{\nu/\alpha_0})\right)\d x \d \mbs'' \d \mbs'.
\end{align*}
Denote
\begin{equation} \label{h:def}
h(x) := x^{\alpha_0} \E [R^{1-\nu_0/\nu} \I(R > x^{\nu/\nu_0})], \quad
x >0.
\end{equation}
Note the limit
\begin{eqnarray} \label{hlim}
\lim_{x \to \infty} h(x)
&=&c_R \, \alpha\lim_{x \to  \infty} x^{\alpha_0}\int_{x^{ \nu/\nu_0}}^\infty \frac{r^{1- \nu_0/\nu}}{r^{1+ \alpha}}  \d r  \
=
\frac{c_R \alpha \nu}{\alpha_0 \nu_0}
=:
\ h_\infty.
\end{eqnarray}
This is a consequence of the tail behavior of $R$ in
\eqref{falpha} with $\alpha<1+\frac{\nu_0}{\nu}$ and $(\alpha -1)\nu = (\alpha_0 -1) \nu_0$. Recalling that $\E(R^{1-\nu_0/\nu})\le \E(R)^{1-\nu_0/\nu}<+\infty$, it follows that $h$ is bounded.
Conditioning by $g^0(\mbs'')$ and using the independence with $R$ inside the expectation we obtain
\begin{align*}
&\tilde j_{0,\lambda}(\theta)\\
&=\lambda^{\nu_0}
\int_{\R^{\nu_0}} \i \theta \phi(\mbs') \int_{\R^{\nu - \nu_0}}  \int_{\R^+}
\Big(\e^{\i \theta\phi(\mbs') x} - 1\Big)
 \E\left(g^0(\mbs'')^{\alpha_0} h\Big(\frac{x \lambda^{\nu_0/\alpha_0}}{g^0(\mbs'')}\Big)  \right)\frac{\d x}{ (x \lambda^{\nu_0/\alpha_0})^{\alpha_0}}\d \mbs'' \d \mbs'\\
&=
\int_{\R^{\nu_0}} \i \theta \phi(\mbs') \d \mbs'  \int_{\R^+}
\Big(\e^{\i \theta\phi(\mbs') x} - 1\Big) \frac{ \d x}{x^{\alpha_0}}
\int_{\R^{\nu - \nu_0}} \E\left(g^0(\mbs'')^{\alpha_0} h\Big(\frac{x \lambda^{\nu_0/\alpha_0}}{g^0(\mbs'')}\Big)\right)  \d \mbs''.
\end{align*}
Then, setting
\begin{equation}\label{c0Xi}
 c_{0,\Xi}:=h_\infty \int_{\R^{\nu - \nu_0}}  \E\left(g^0(\mbs'')^{\alpha_0}  \right)\d \mbs'',
\end{equation}
we infer from the dominated convergence theorem that
\begin{eqnarray}\label{j00lambda}
\tilde j_{0,\lambda}(\theta)
&\to&j_0(\theta) := c_{0,\Xi}
\int_{\R^{\nu_0}} \i \theta \phi(\mbs') \d \mbs'  \int_{\R^+}
\Big(\e^{\i \theta\phi(\mbs') x} - 1\Big) \frac{ \d x}{x^{\alpha_0}},
\end{eqnarray}
$j_0(\theta)$ being the log-characteristic function of the $\alpha_0$-stable random variable
$L_{\alpha_0}(\phi)$ in \eqref{limvolX0}.

Let us prove that
$$
\lim_{\lambda \to \infty}  |j_{0,\lambda}(\theta) -  \tilde j_{0,\lambda}(\theta)|=0.
$$
We follow \eqref{jnew} and  \eqref{jnew1} in the proof of Theorem~\ref{thm1}.
For simplicity we give the proof assuming $\phi \in \Phi_0$ uniformly continuous on $\R^{\nu_0}$ to avoid
the approximation step using Lusin's theorem.
Analogously,
introduce
\begin{eqnarray*}
\xi_{0,\lambda} (\mbs',\mbs'')
&:=&\lambda^{-\nu_0/\alpha_0}
\int_{\R^{\nu_0}}\phi\big(\frac{\mbt'}{\lambda} +  \mbs'\big) \I ((\mbt',\boldsymbol{0}) \in (\boldsymbol{0},\mbs'')  + \Xi) \d \mbt', \\
\tilde \xi_{0,\lambda} (\mbs',\mbs'') &:=&  \lambda^{-\nu_0/\alpha_0} \phi(\mbs') \Leb_{\nu_0}(\Xi_{\mbs''}), \\
\eta_{0,\lambda} (\mbs',\mbs'') &:=& \xi_{0,\lambda} (\mbs',\mbs'') - \tilde \xi_{0,\lambda} (\mbs',\mbs'') \\
&=& \lambda^{-\nu_0/\alpha_0}
\int_{\R^{\nu_0}}\big(\phi\big(\frac{\mbt'}{\lambda} +  \mbs'\big) -  \phi(\mbs')\big)
\I ((\mbt',\boldsymbol{0}) \in (\boldsymbol{0},\mbs'') + \Xi) \d \mbt'.
\nn
\end{eqnarray*}
Without loss of generality we can assume that $\phi(\mbs')>0$. We have, integrating by parts,
\begin{eqnarray*}
j_{0,\lambda}(\theta)
&=&\lambda^{\nu_0} \int_{\R^{\nu_0} \times \R^{\nu - \nu_0}}
\E \Psi \big(\theta \xi_{0,\lambda} (\mbs',\mbs'') \big)
\d \mbs' \d \mbs''  \\
&=&\i \theta  \int_{\R^{\nu_0}} \d \mbs' \int_0^\infty (\e^{\i \theta x} -1) \d x \  \lambda^{\nu_0}  \int_{\R^{\nu - \nu_0}}
\P(\xi_{0,\lambda} (\mbs',\mbs'') >x ) \d \mbs''.
\end{eqnarray*}
Following \eqref{limla}, let us check that for all $x>0$,
\begin{eqnarray}\label{limla00}
&&\lim_{\lambda \to \infty} \lambda^{\nu_0} x^{\alpha_0} \int_{\R^{\nu-\nu_0}} \P(\xi_{0,\lambda} (\mbs',\mbs'') > x) \d \mbs''
\ = \
\tilde c_0 \phi(\mbs')^{\alpha_0}
\end{eqnarray}
which can be compared to \eqref{j00lambda} by noting that
\begin{eqnarray*}
&&\phi(\mbs') \int_0^\infty
\big(\e^{\i \theta\phi(\mbs') x} - 1\big) \frac{ \d x}{x^{\alpha_0}}\\
&&= |\phi(\mbs')|^\alpha \Big(\int_0^\infty
\big(\e^{\i \theta x} - 1\big) \frac{ \d x}{x^{\alpha_0}} \I(\phi(\mbs') >0) +\int_0^\infty
\big(\e^{-\i \theta x} - 1\big) \frac{ \d x}{x^{\alpha_0}} \I(\phi(\mbs') <0)\Big).
\end{eqnarray*}

To show  \eqref{limla00},
similarly as in \eqref{xieta}, for any fixed $\gamma >0$ we evaluate
the integral on the left hand side as
\[
\int_{\R^{\nu-\nu_0}} \P\big(\xi_{0,\lambda} (\mbs',\mbs'') > x\big) \d \mbs''
\le\int_{\R^{\nu-\nu_0}}\Big(\P\big(\tilde \xi_{0,\lambda} (\mbs',\mbs'') > (1-\gamma) x\big)
+ \P\big(|\eta_{0,\lambda}(\mbs',\mbs'')|  >\gamma x\big)\Big) \d \mbs'',
\]
\[
\int_{\R^{\nu-\nu_0}} \P\big(\xi_{0,\lambda} (\mbs', \mbs'') > x\big) \d \mbs''
\ge\int_{\R^{\nu-\nu_0}} \Big(\P\big(\tilde \xi_{0,\lambda} (\mbs',\mbs'') > (1+\gamma) x\big)
- \P\big(|\eta_{0,\lambda}(\mbs',\mbs'')|  >\gamma x\big)\Big) \d \mbs''.
\]
With \eqref{j00lambda} in mind and taking $\gamma >0$ arbitrary small, this
reduces the proof of \eqref{limla00} to
\begin{eqnarray}\label{limetala0}
&&\lim_{\lambda  \to \infty} \lambda^{\nu_0} x^{\alpha_0} \int_{\R^{\nu-\nu_0}}
\P\big(|\eta_{0,\lambda} (\mbs',\mbs'')| > x\big)\d \mbs'' \  = \ 0.
\end{eqnarray}
Proceeding similarly to \eqref{limetala}, by
uniform continuity of $\phi$, for any $\epsilon >0$ there is a $\tau >0$ such  that
$$
\sup_{\mbs' \in B_K} \sup_{|\mbt'| \le \tau \lambda}  \big|
\phi\big(\frac{\mbt'}{\lambda} +  \mbs'\big) -  \phi(\mbs')\big| < \epsilon
$$
uniformly in  $\lambda >0$.
Therefore,
\begin{eqnarray*}
|\eta_{0,\lambda}(\mbs', \mbs'')|
&\le&\lambda^{-\nu_0/\alpha_0}
\int_{\R^{\nu_0}} |\phi(\frac{\mbt'}{\lambda} + \mbs') -
\phi(\mbs')|  \I ((\mbt',\boldsymbol{0}) \in ({\boldsymbol{0}},\mbs'')  + \Xi)
\d \mbt' \\
&\le&\lambda^{-\nu_0/\alpha_0} ( \epsilon \Leb_{\nu_0} (\Xi_{\mbs''}) +
C \Leb_{\nu_0} (\Xi_{\mbs''} \cap \{|\mbt' | > \tau \lambda \})).
\end{eqnarray*}
Consider the first term on the right hand side above.
Since
$$\Leb_{\nu_0}(\Xi_{\mbs''})
= R^{\nu_0/\nu} g^0(\mbs''/R^{1/\nu})
\le  C R^{\nu_0/\nu}\I(|\mbs''| < R^{1/\nu}),$$
using the boundedness $h(x) \le C $, see \eqref{hlim}, we get that
\begin{eqnarray*}
 \lambda^{\nu_0} x^{\alpha_0}\int_{\R^{\nu - \nu_0}}
\P\big(\lambda^{-\nu_0/\alpha_0} \epsilon \Leb_{\nu_0} (\Xi_{\mbs''}) > x\big)\d \mbs''
&\le&C \lambda^{\nu_0} x^{\alpha_0} \E [R^{1 - \nu_0/\nu} \I(R > \lambda^{\nu/\alpha_0} (x/\epsilon)^{\nu/\nu_0})] \\
&\le&C \epsilon^{\alpha_0}
\end{eqnarray*}
that can be made arbitrarily small with $\epsilon >0$
uniformly in $x$ and $\lambda$.

Next, consider
\begin{eqnarray*}
\Leb_{\nu_0} (\Xi_{\mbs''} \cap \{|\mbt' | > \tau \lambda \})
&=&
\Leb_{\nu_0} ( R^{1/\nu} \Xi^0 \cap \{\mbt'' = \mbs''\}  \cap \{|\mbt' | > \tau \lambda \}) \\
&=&
R^{\nu_0/\nu} \Leb_{\nu_0} (\Xi^0 \cap \{\mbt'' = \mbs''/R^{1/\nu}\}  \cap \{|\mbt' | > \tau \lambda/R^{1/\nu} \})\\
&\le&
\begin{cases}
R^{\nu_0/\nu} g^0(\mbs''/R^{1/\nu}) &\text{ if }\tau \lambda < R^{1/\nu}, \\
0  &\text{ if }\tau \lambda > R^{1/\nu}.
\end{cases}
\end{eqnarray*}
Therefore, using that $g^0$ is bounded with bounded support and $h$ given by \eqref{h:def} is bounded,
for any $\tau >0$,
\begin{eqnarray*}
&&\int_{\R^{\nu-\nu_0}}
\P\big(\Leb_{\nu_0} (\Xi_{\mbs''} \cap \{|\mbt' | > \tau \lambda \})  > x\lambda^{\nu_0/\alpha_0}/C\big)\d \mbs'' \\
&&\le\ \E \int_{\R^{\nu-\nu_0}}\I(\tau \lambda < R^{1/\nu})
 \I \big(R^{\nu_0/\nu} g^0(\mbs''/R^{1/\nu})  > x\lambda^{\nu_0/\alpha_0}/C\big)\d \mbs'' \\
&&\le  C  \left((x^{\nu/\nu_0}\lambda^{\nu/\alpha_0}) \vee (\tau \lambda)^\nu )\right)^{-\alpha_0} h\left((x^{\nu/\nu_0}\lambda^{\nu/\alpha_0}) \vee (\tau \lambda)^\nu \right)\\
&&\le  C   \left((x^{\nu/\nu_0}\lambda^{\nu/\alpha_0}) \vee (\tau \lambda)^\nu \right)^{-\alpha_0},
\end{eqnarray*}
implying
$\lim_{\lambda \to \infty}\lambda^{\nu_0} x^{\alpha_0} \int_{\R^{\nu-\nu_0}}
\P\big(\Leb_{\nu_0} (\Xi_{\mbs''} \cap \{|\mbt' | > \tau \lambda \})  > x\lambda^{\nu_0/\alpha_0}/C\big)\d \mbs''
= 0$
and ending the proof of \eqref{limetala0}, and thus of \eqref{limvolX0}.

\medskip

\noi Let us turn to Step 2, that is the proof of the fdd convergence \eqref{limvolXk}. Following the proof of Step 2 in Theorem~\ref{thm1},
it suffices to check
\begin{equation}\label{relZ20}
{\rm Var}({\mathcal Z}_{0,\lambda}(\phi))
=\int_{\R^{2\nu_0}} \phi(\mbt_1'/\lambda) \phi(\mbt_2'/\lambda)
{\rm Cov}({\mathcal Z}(\mbt_1', \boldsymbol{0}),  {\mathcal Z}(\mbt_2', \boldsymbol{0})) \d \mbt_1' \d \mbt_2' \ = \
o(\lambda^{2\nu_0/\alpha_0}),
\end{equation}
where ${\mathcal Z}(\mbt) \equiv {\mathcal Z}_k(\mbt)$, $\mbt \in \R^\nu $ is the same as in \eqref{relZ2}. Relation \eqref{relZ20} follows similarly to \eqref{relZ2}, using \eqref{relZ1}, \eqref{covXRG}, and the fact that for $\beta := 2\nu (\alpha -1)$,
\[
2 \nu_0 - \beta=2\nu_0(2-\alpha_0)<2\nu_0/\alpha_0
\]
since for $1<\alpha<1+\frac{\nu_0}{\nu}$, $\alpha_0\in(1,2)$. See also \cite[Proposition 5, (56)]{pils2024}.

\smallskip

\noi (ii) The proof is similar to Theorem~\ref{thm2}. Essentially, we need to check only
$\int_{\R^{\nu_0}} r_X(\mbt',\boldsymbol{0}) \d \mbt' < \infty $. This is immediate from \eqref{covXRG} and the boundedness of $r_X$, yielding
\begin{eqnarray*}
\int_{\R^{\nu_0}} r_X(\mbt',\boldsymbol{0})\d \mbt'&\le&C\int_{\R^{\nu_0}}
(1 \wedge |\mbt'|)^{-\nu(\alpha -1)} \d \mbt' < \infty
\end{eqnarray*}
for $\nu (\alpha - 1) > \nu_0 $, or
$\alpha > 1 + \frac{\nu_0}{\nu}$. \hfill $\Box$

\section{Numerical illustrations}\label{sec:num}
\begin{figure}[h]
\begin{center}
\begin{tabular}{lll}
\hspace{-0.5cm}\includegraphics[width=5.5cm]{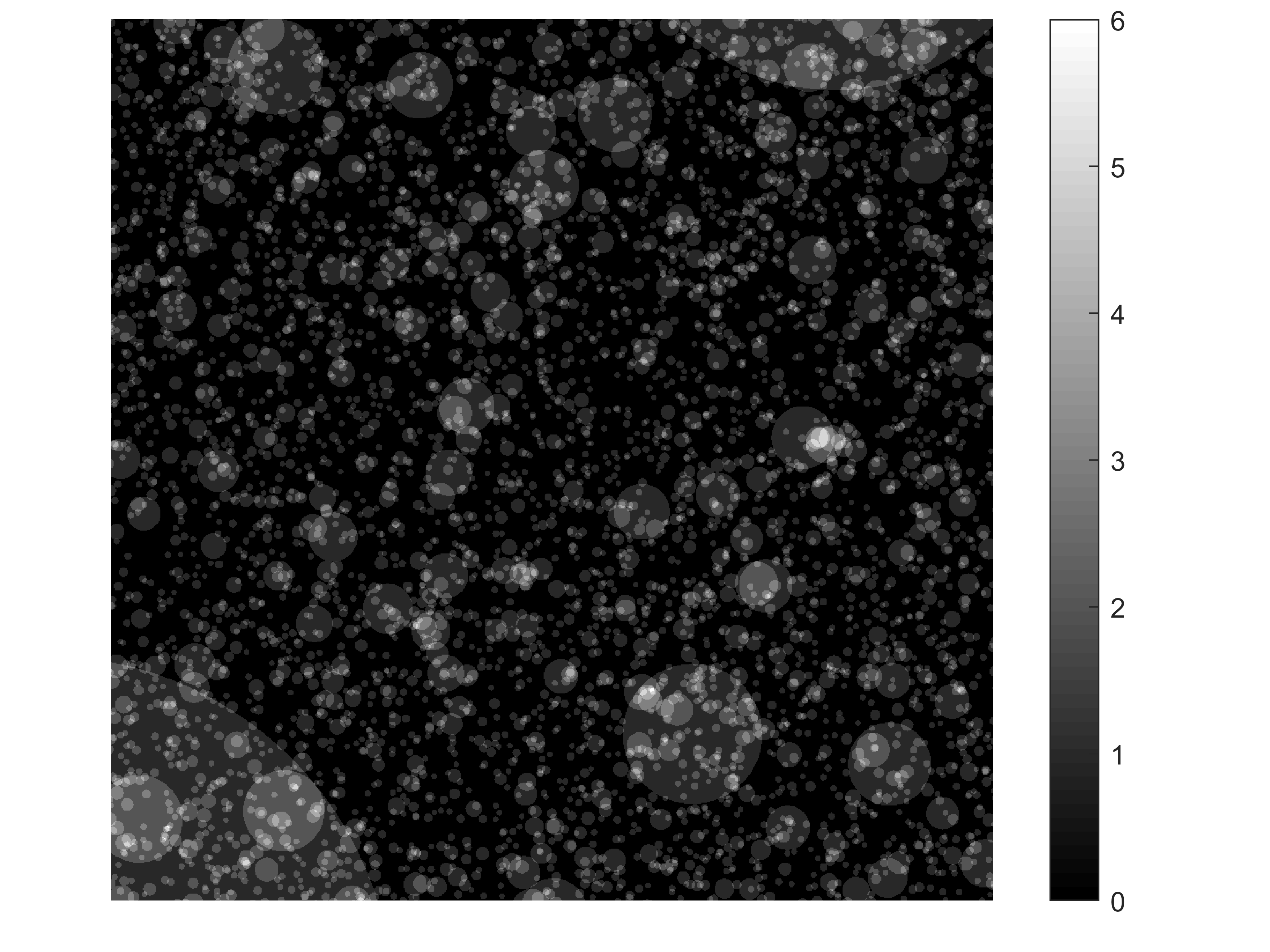} &
\includegraphics[width=5.5cm]{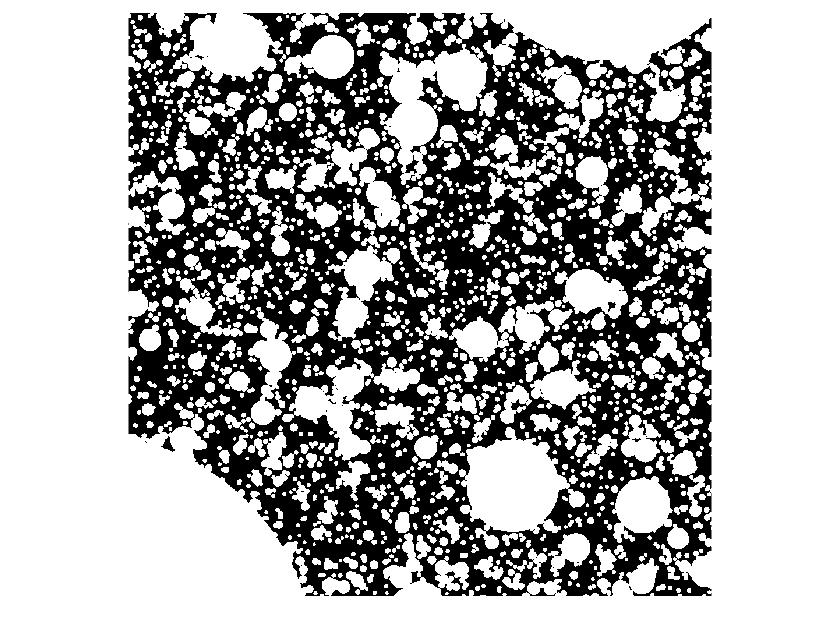} &
\includegraphics[width=5.5cm]{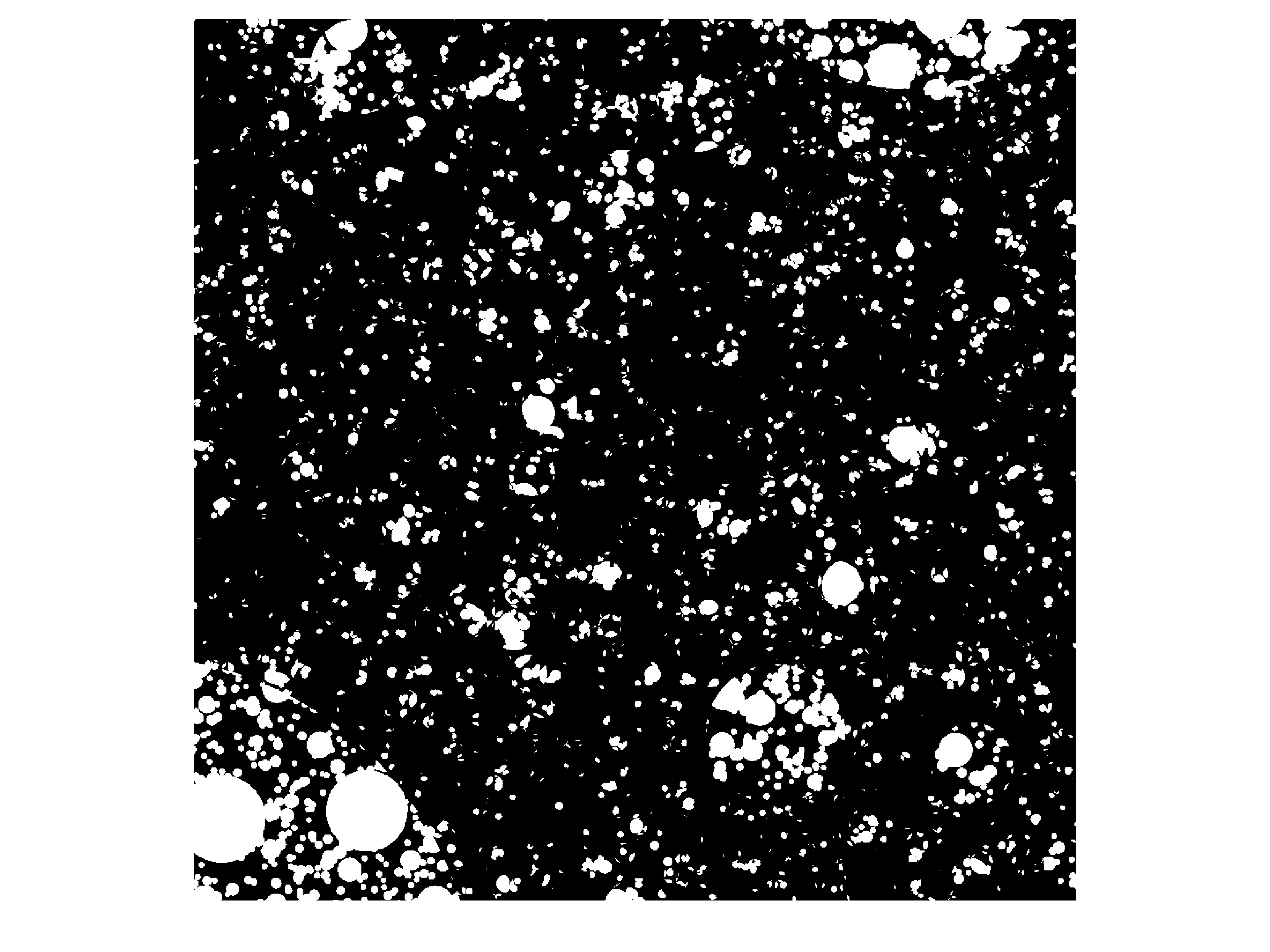}
\end{tabular}
\end{center}
\caption{\small Random grain model for $\alpha=1.3$, $a_p=0.005$ and $c=10a_p$. Left: a sample of the random field with color-bar for values; middle: the associated Boolean field; right: the excursion set for $k=2$. Pixel values equal to one are drawn in white.}
\label{SN1:img}
\end{figure}

In this section we consider simulation of the random homothetic grain random field in dimension $\nu=2$ where $\Xi=R^{1/\nu}B$, with
$B=B(0,1)$ is the Euclidean ball and $R$ is a random variable given by
$$R=c U^{-1/\alpha},$$
for some $c>0$, $\alpha\in (1,2)$ and $U$ a uniform random variable on $(0,1)$. Note that for all $x>c$ one has $\P(R>x)=c^\alpha x^{-\alpha}$ such that $c_R=c^\alpha$ and $c_\Xi=(c\pi)^\alpha$ according to Proposition \ref{corRGhomo}. Note also that the mean volume of the ball is given by
\begin{equation*}
    \mu=c\pi\E(R)  =c\pi\frac{\alpha}{\alpha-1}.
\end{equation*} To illustrate our results we fix a rectangular observation window $\lambda A$ consisting of $N \times N$ pixels, each with volume $a_p$, where
$N = 1000$ and $a_p = 0.005$. Hence
$\Leb_\nu (\lambda A) = N^2 a_p = 5000, \lambda = \sqrt{5000}$
for $A = [0,1]^2$. Then we approximate the Lebesgue measure of 
a subset of $\lambda A$ by counting the number of pixels it contains and multiplying by $a_p$. Figure \ref{SN1:img} shows a simulation of the random grain model described above and its excursion sets for $k=1$, which corresponds to the Boolean set ${\mathcal X}$, and $k=2$.

\begin{figure}[h]
\begin{center}
\includegraphics[width=7cm]{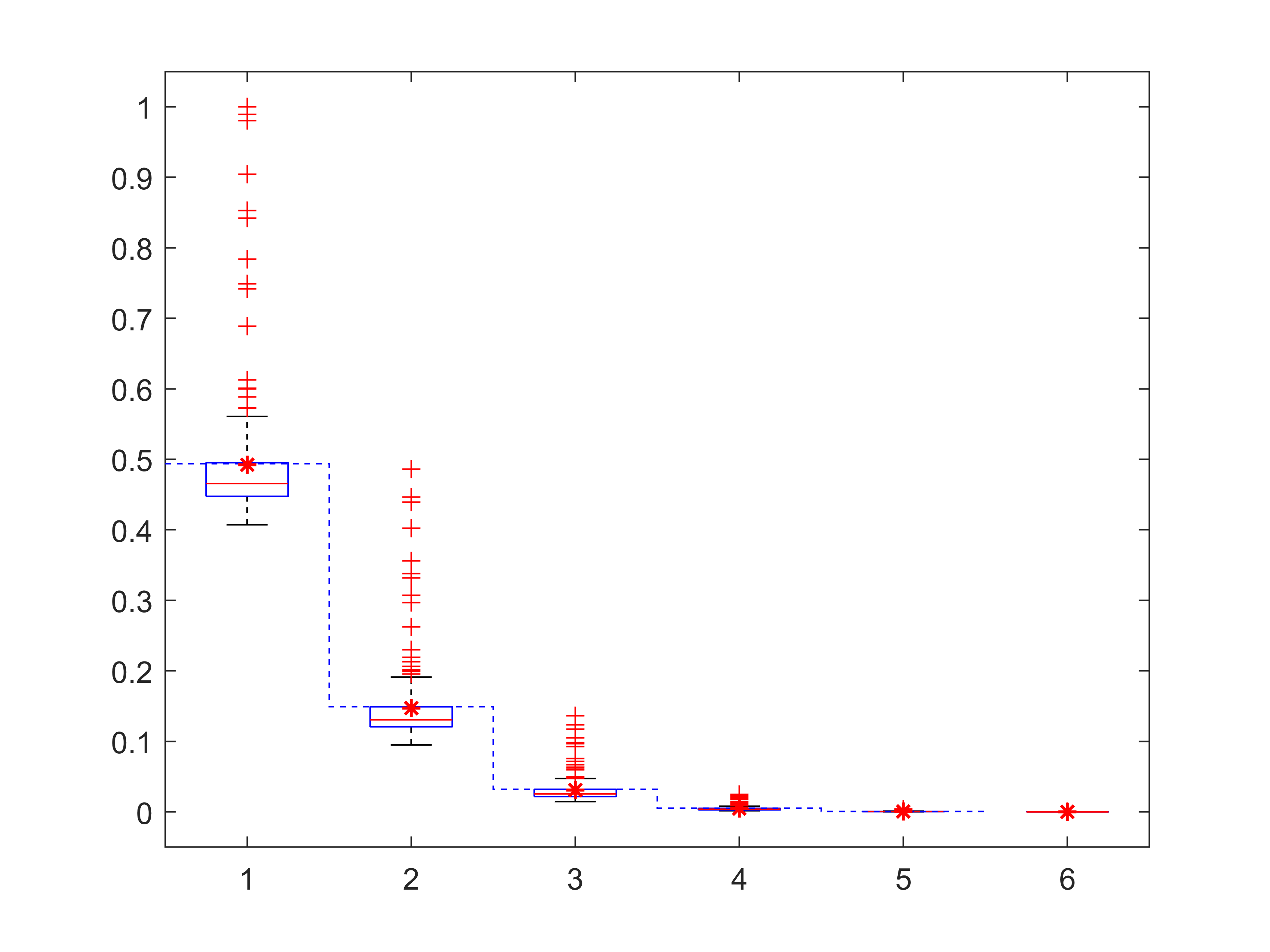}
\end{center}
\caption{\small Random grain model for $\alpha=1.3$, $a_p=0.005$ and $c=10a_p$. Boxplots of the estimated volume fraction of the excursion sets of the random grain  model, based on 200 independent realizations, $\{\widehat{p}_{\lambda,k}; 1\le k\le 6\}$ given by \eqref{volumefe}. The red stars indicate the empirical mean value. The dotted blue stairs represent the theoretical values $\{p_k; 1\le k\le 6\}$.  }
\label{SN1:Boxplot}
\end{figure}

In view of \eqref{volumef} and $\Leb_\nu (A)=1$, we take
\begin{equation}\label{volumefe}
\widehat p_{\lambda, k} :=\lambda^{-\nu} \widehat X_{\lambda,k} (A)= \lambda^{-\nu}\Leb_\nu (\{\mbt \in \lambda A; X(\mbt)\ge k\}),
\end{equation}
as the estimator of the volume fraction $p_k$ of the corresponding excursion set that is given by
$$
p_k=\P(X(\mbt)\ge k)=1-e^{-\mu}\sum_{j=0}^{k-1}\frac{\mu^j}{j!},$$
We refer to Figure \ref{SN1:Boxplot} for a comparison between empirical and theoretical values.

\begin{figure}[h]
\begin{center}
\begin{tabular}{ccc}
\hspace{-0.5cm}\includegraphics[width=5.5cm]{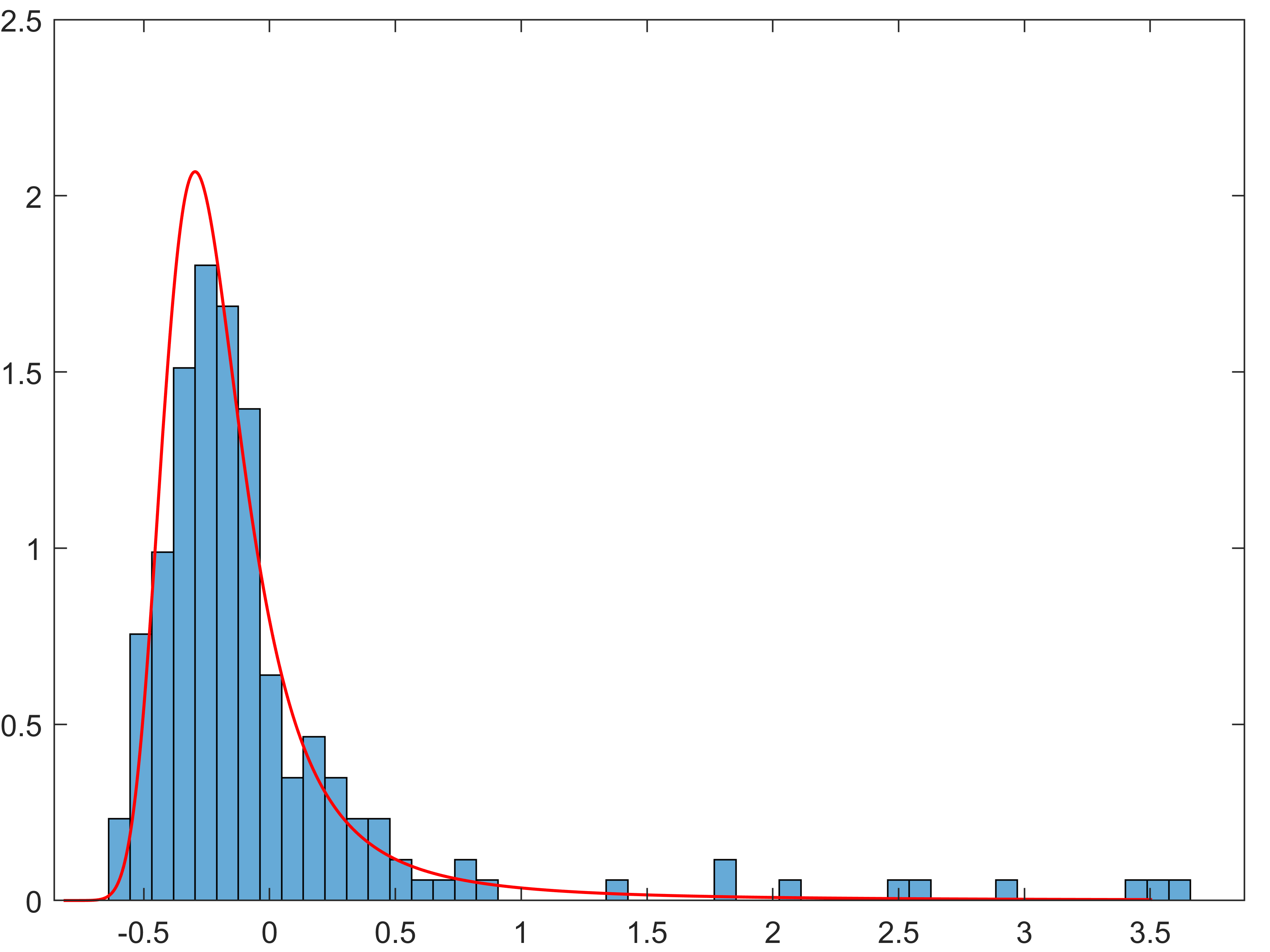} &
\includegraphics[width=5.5cm]{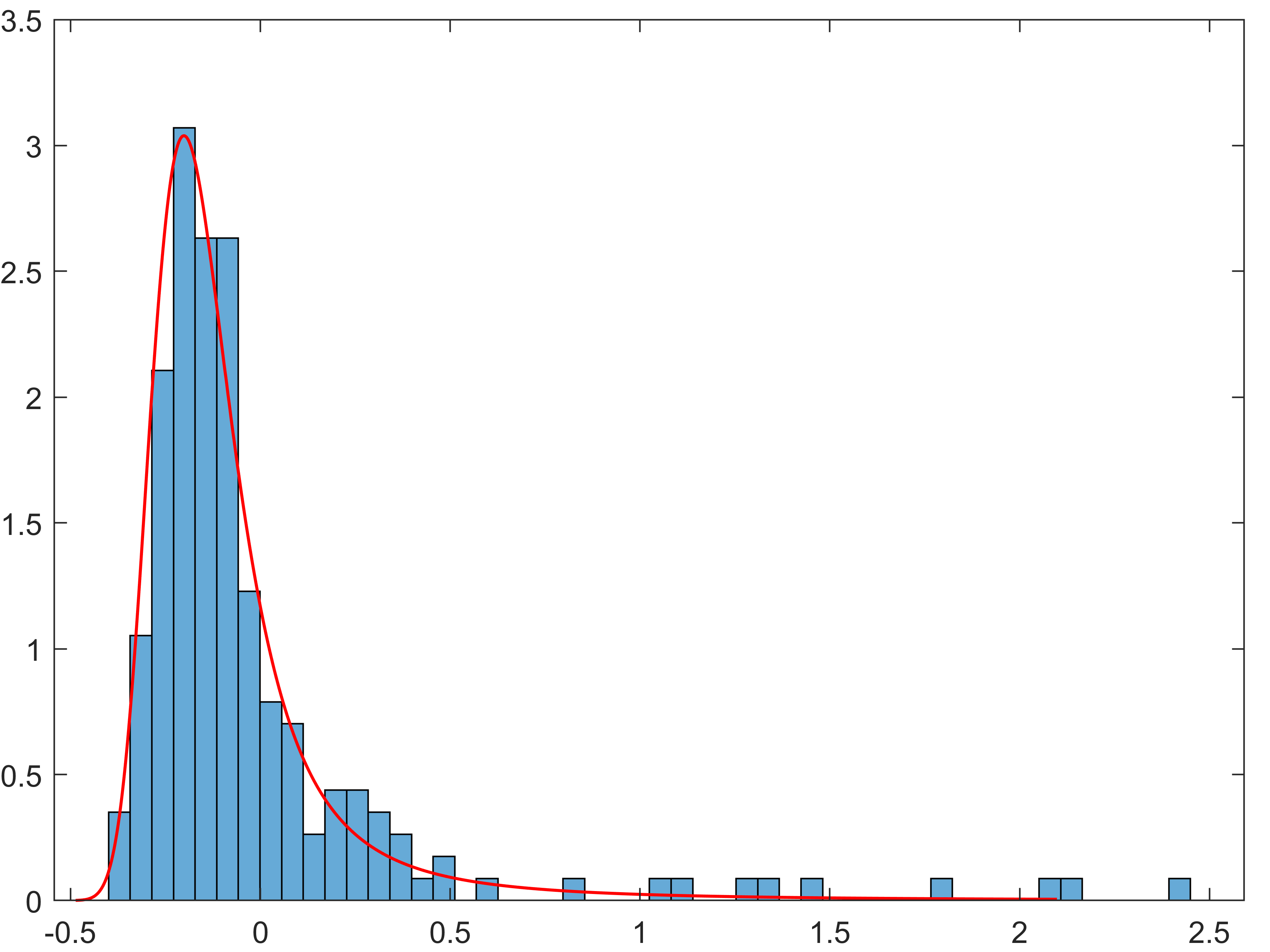} &
\includegraphics[width=5.5cm]{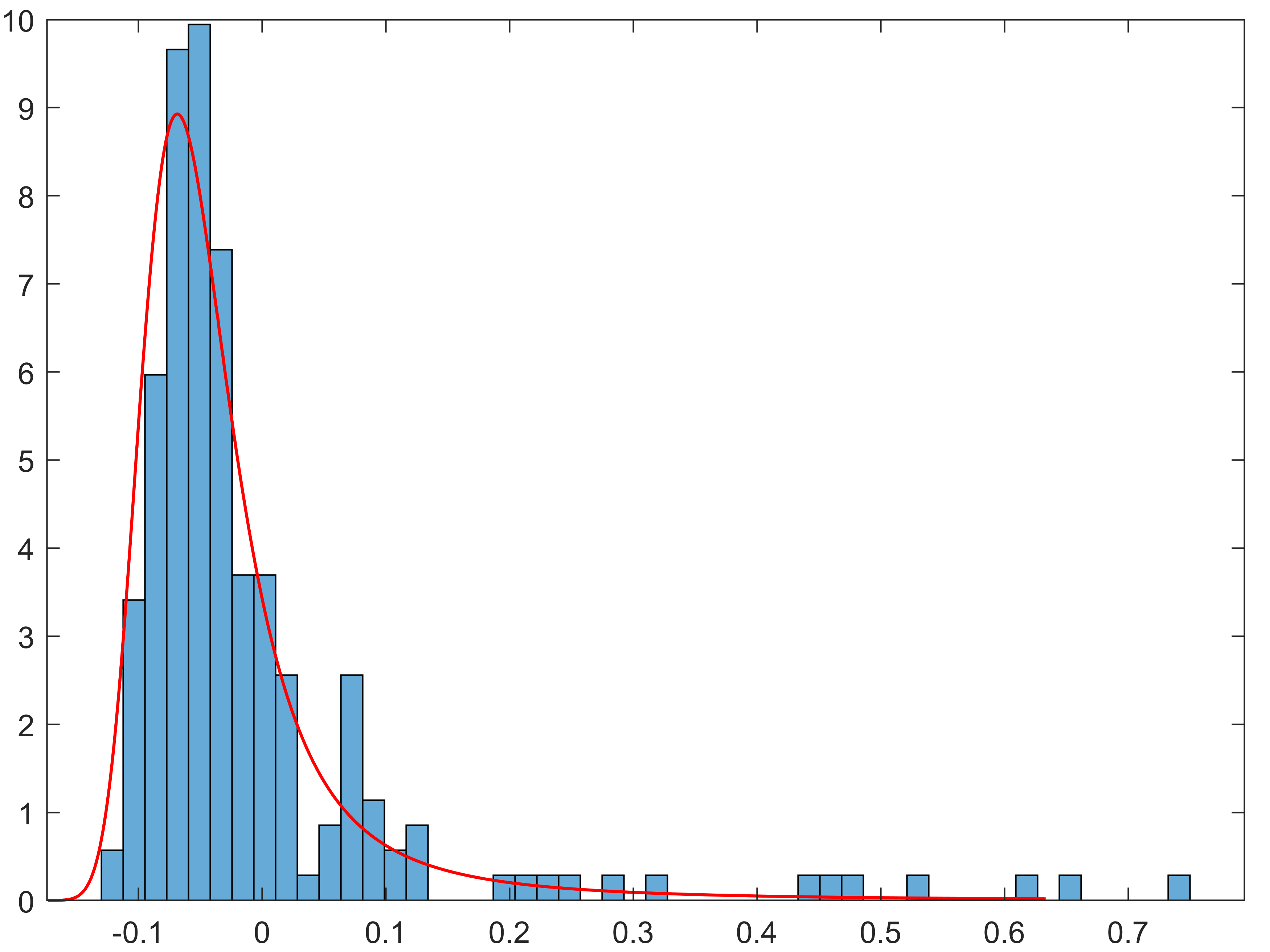}\\
$k=1$&$k=2$&$k=3$
\end{tabular}
\end{center}
\caption{\small Random grain model for $\alpha=1.3$, $a_p=0.005$ and $c=10a_p$. Histograms of the statistic $\lambda^{\nu-\nu/\alpha}\lambda^{\nu-\nu/\alpha} (\widehat p_{\lambda,k}-p_k)$ for $k\in\{1,2,3\}$, based on $200$ independent realizations. The red curves show the theoretical asymptotic  probability distributions in view of \eqref{volumefe1}.}
\label{SN1:Histo}
\end{figure}

Then, from Theorem \ref{thm1} and Corollary \ref{Cor:thm1}
\begin{equation}\label{volumefe1}
\lambda^{\nu-\nu/\alpha} (\widehat p_{\lambda,k}-p_k)=\lambda^{-\nu/\alpha}\left( \widehat X_{\lambda,k} (A)-  \E \widehat X_{\lambda,k} (A)\right)\overset{\d}{\longrightarrow}\e^{-\mu} \frac{\mu^{k-1}}{(k-1)!} L_\alpha (A),
\end{equation}
where,
 recalling \eqref{cf:cor1},
$$
\E \e^{\i \theta L_\alpha(A)} =
\exp\Big\{ \i c_\Xi \theta \Leb_\nu(A) \int_{\R_+}
\big(\e^{\i \theta x} -1 \big)x^{-\alpha} \d x\Big\} .
$$
According to  (3.9) in chapter XVII of \cite{feller1991introduction}, for any $a \in (0,1)$, and $\theta>0$
$$ \int_{\R_+}
\big(\e^{\i \theta x} -1 \big)x^{-a-1} \d x =\theta^a \frac{\Gamma(2-a)}{a(a-1)}e^{-i\pi a/2}.$$
Hence, taking $\alpha=1+a$,
one has
$$i\theta \int_{\R_+}
\big(\e^{\i \theta x} -1 \big)x^{-\alpha} \d x=\theta^\alpha\frac{\Gamma(2-\alpha)}{\alpha-1}e^{-i\pi\alpha/2},$$
from which we deduce that, for all $\theta\in\R$,
$$i\theta \int_{\R_+}
\big(\e^{\i \theta x} -1 \big)x^{-\alpha} \d x=|\theta|^\alpha\frac{\Gamma(2-\alpha)}{\alpha-1}\cos(\pi\alpha/2)\left(1-i\mbox{sgn}(\theta)\tan(\pi\alpha/2)\right).$$
Therefore, in view of Definition 1.1.6 of \cite{samorodnitsky1994stable},
the $\alpha$ stable random variable $L_\alpha(A)$ follows a stable distribution of parameters $\alpha$, $\beta=1$, $\delta=0$ and scale parameter
\begin{equation*}
\gamma=\left(c_\Xi \Leb_\nu(A) \frac{\Gamma(2-\alpha)}{\alpha-1}\left|\cos(\pi\alpha/2)\right|\right)^{1/\alpha}.
\end{equation*}
Using the Matlab package STBL \cite{veillette_stbl} this will allow us to  compare empirical results with the theoretical probability distribution in Figure \ref{SN1:Histo}.

Next, we consider empirical estimates of the volume fractions defined above along lines parallel to the coordinate axes. We have $\nu_0=1$, $\alpha_0=2\alpha-1$ and we define
\begin{equation}\label{volumeLfe}
L\widehat p_{\lambda, k} :=\lambda^{-\nu_0} \widehat X_{\lambda,k} (A\cap H_{\nu_0})= \lambda^{-\nu_0}\Leb_{\nu_0} (\{\mbt \in \lambda A\cap H_{\nu_0}; X(\mbt)\ge k\}),
\end{equation} 
as the estimator of the volume fraction computed from a single  line of the image.
\begin{figure}[h]
\begin{center}
\begin{tabular}{cccc}
\hspace{-0.5cm}\includegraphics[width=5cm]{
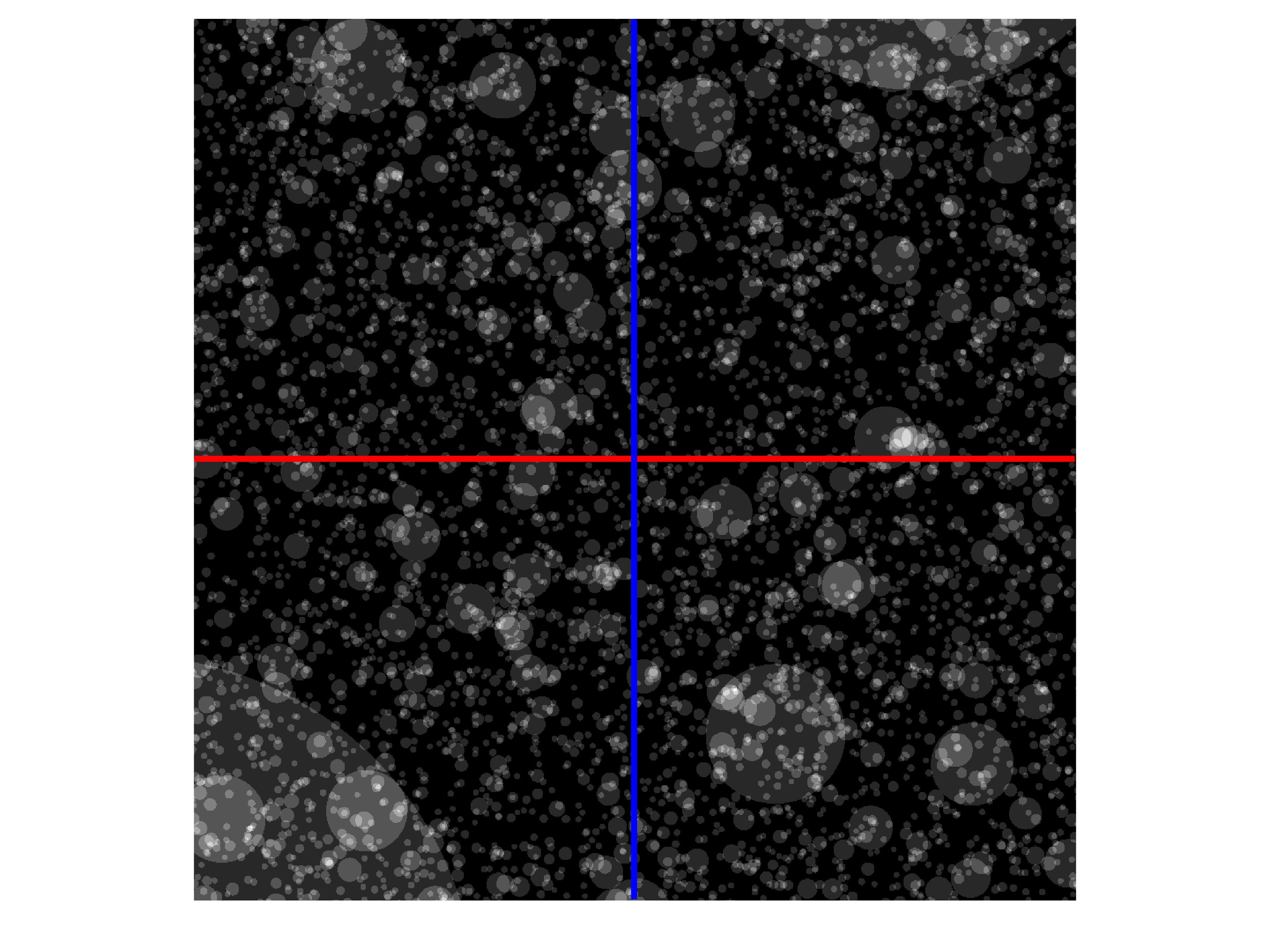} &
\includegraphics[width=5cm]{
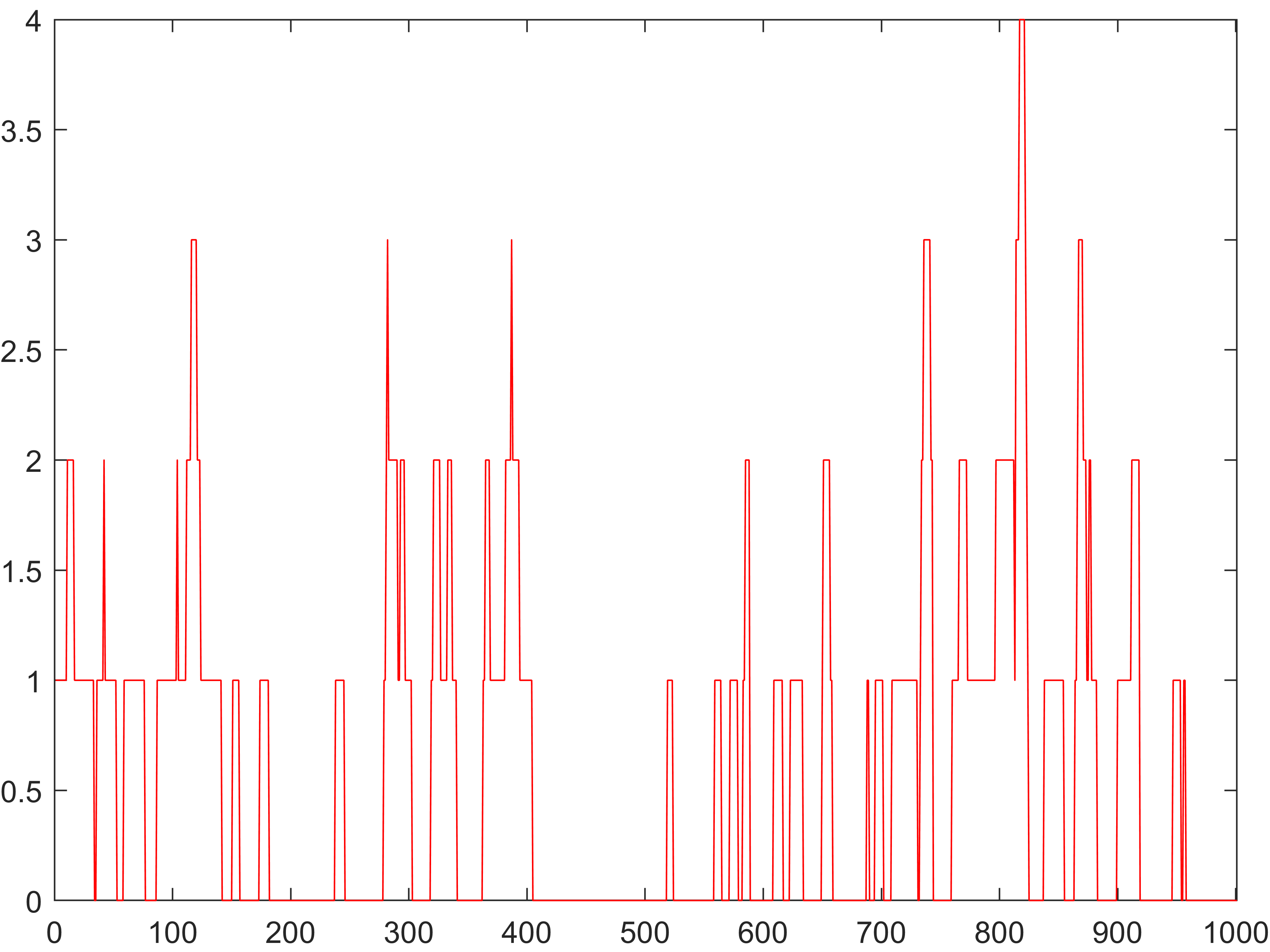} &&
\includegraphics[width=5cm]{
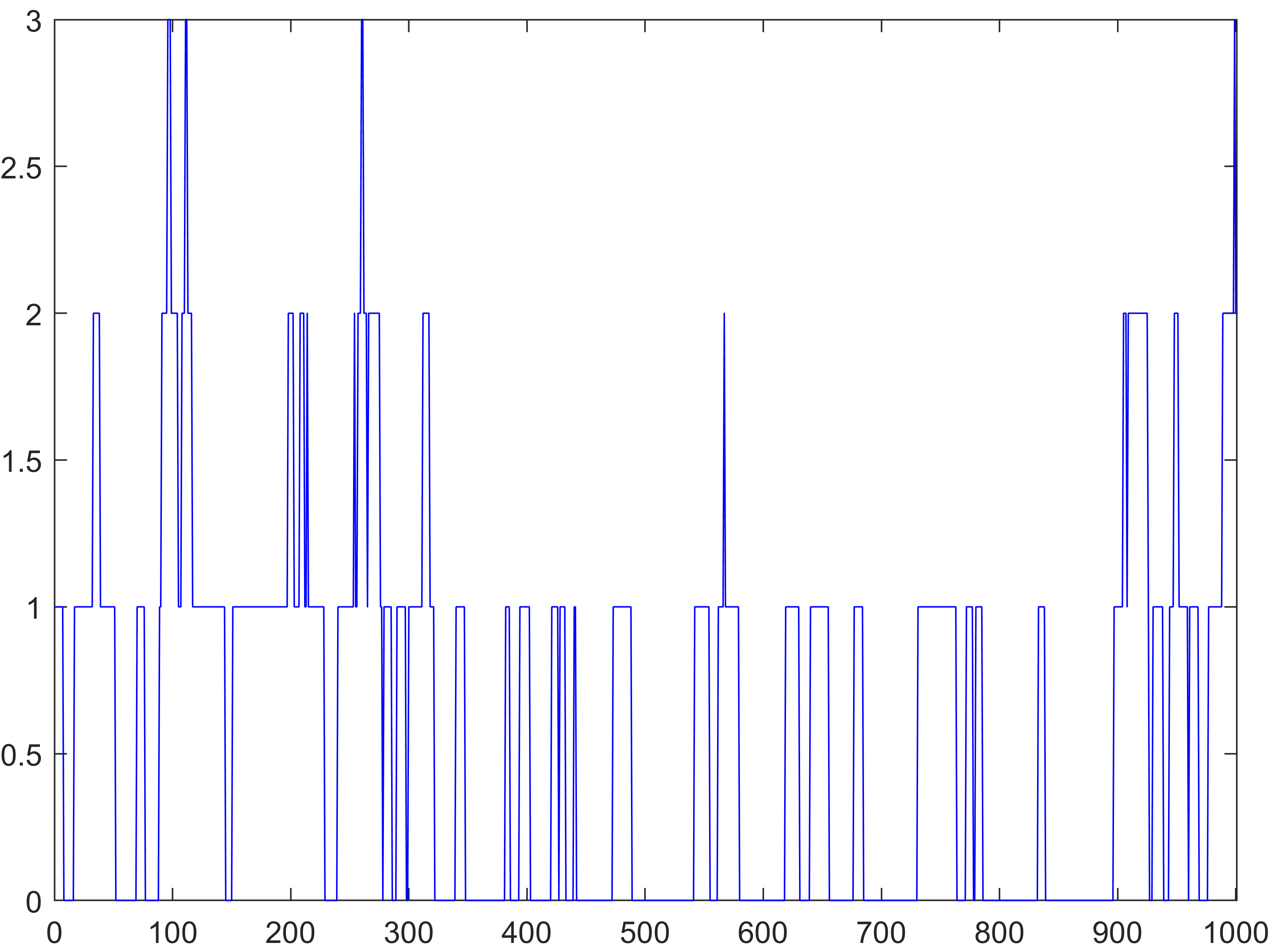}\\
&horizontal&&vertical
\end{tabular}
\end{center}
\caption{\small Left: random grain model for $\alpha=1.3$, $a_p=0.005$ and $c=10a_p$. Middle: Observations of the random field along the horizontal line. Right: Observations of the random field along the vertical line.}
\label{fig:extraction}
\end{figure}

\begin{figure}[h]
\begin{center}
\begin{tabular}{cc}
\includegraphics[width=7cm]{
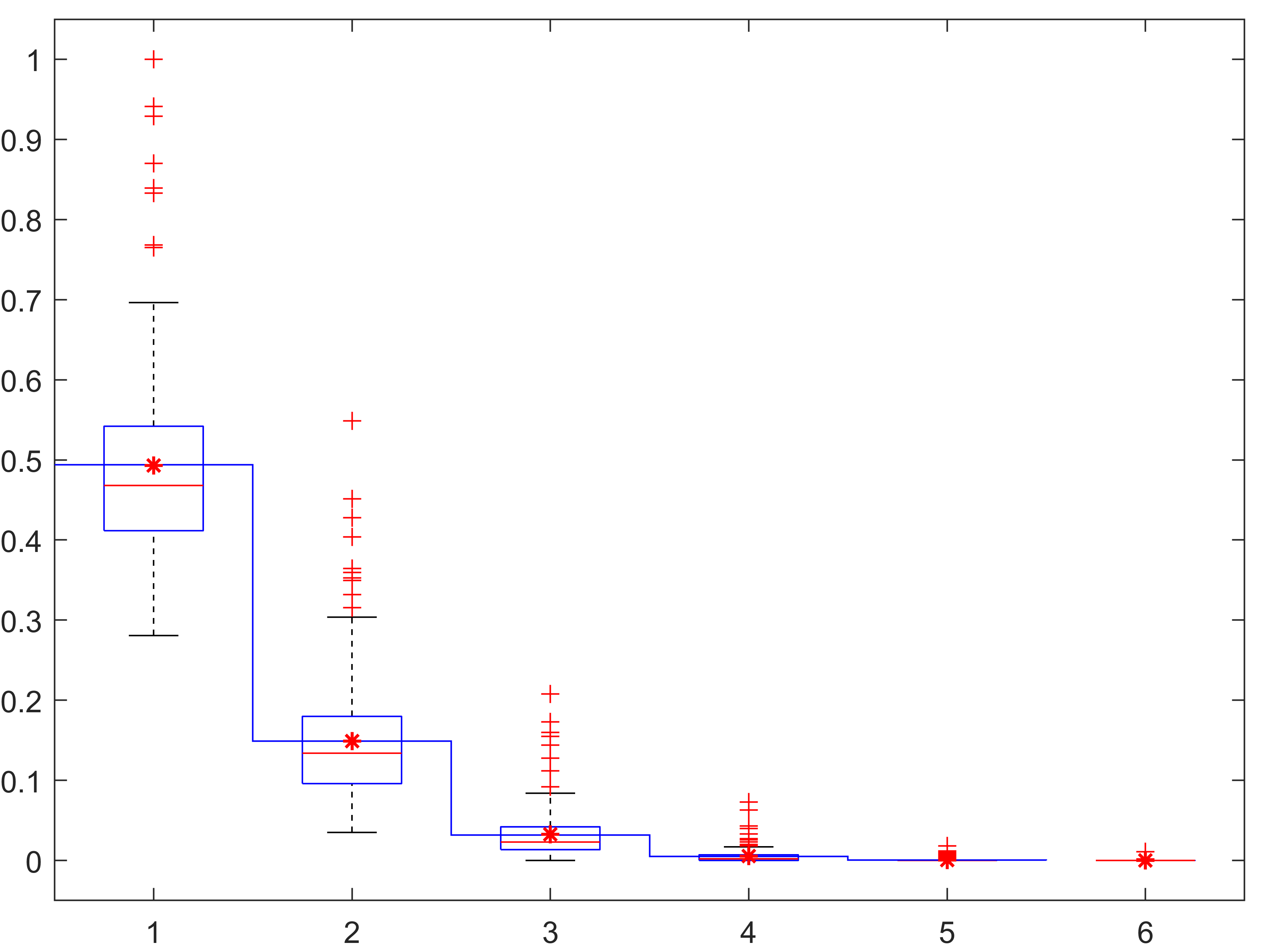} &
\includegraphics[width=7cm]{
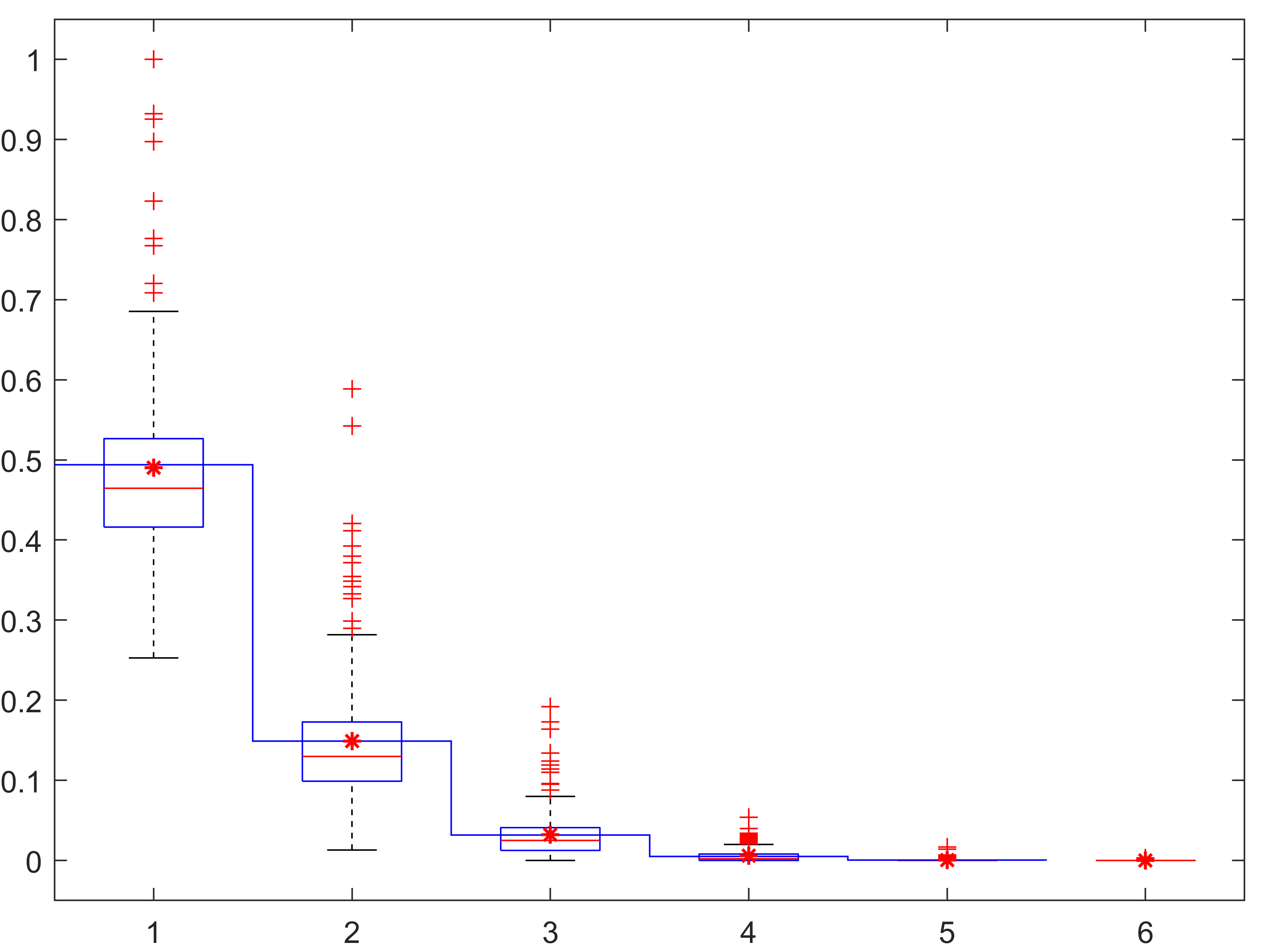}\\
horizontal&vertical
\end{tabular}
\end{center}
\caption{\small Random grain model for $\alpha=1.3$, $a_p=0.005$ and $c=10a_p$.
Boxplots of the estimated volume fraction of the excursion sets of the restrictions along lines, based on 200 independent realizations,  $\{L\widehat{p}_{\lambda,k}; 1\le k\le 6\}$ given by \eqref{volumeLfe}. The red stars indicate the empirical mean value. The dotted blue stairs represent the theoretical values $\{p_k; 1\le k\le 6\}$. Left: horizontal lines, right: vertical ones.}
\label{LSN1:Boxplot}
\end{figure}

Then, for $\alpha_0<2$, taking horizontal or vertical lines (see Figure \ref{fig:extraction}), such that $\Leb_{\nu_0}(A\cap H_{\nu_0})=1$,  Theorem \ref{thm3} and Corollary \ref{cor3} yield
\begin{eqnarray*}
\lambda^{\nu_0-\nu_0/\alpha_0} (L\widehat p_{\lambda,k}-p_k)&=&\lambda^{-\nu_0/\alpha_0}\left( \widehat X_{0,\lambda,k} (A\cap H_{\nu_0})-  \E \widehat X_{0,\lambda, k} (A\cap H_{\nu_0})\right)\\
&\overset{\d}{\longrightarrow}&\e^{-\mu} \frac{\mu^{k-1}}{(k-1)!} L_{\alpha_0} (A\cap H_{\nu_0}).\nonumber
\end{eqnarray*}
The  $\alpha_0$-stable distribution $L_{\alpha_0}(A)$ has parameters given by
$\alpha_0$, $\beta_0=1$, $\delta_0=0$ and
\begin{equation}\label{scale:par0}
\gamma_0=\left(c_{0,\Xi} \Leb_\nu(A) \frac{\Gamma(2-\alpha_0)}{\alpha_0-1}\left|\cos(\pi\alpha_0/2)\right|\right)^{1/\alpha_0}.
\end{equation}
In view of \eqref{c0Xi} we may compute, using the Beta function,
$$\int_{\R}\E\left(g^0(s)^\alpha_0\right) \d s=2^{\alpha_0}\int_{-1}^1(1-s^2)^{\alpha_0/2}\d s =2^{2\alpha_0+1}B\left(\frac{\alpha_0+2}{2},\frac{\alpha_0+2}{2}\right),$$
and $h_\infty=2\frac{\alpha}{\alpha_0}c^\alpha$ such that we explicitly have
$$c_{0,\Xi}=2^{2\alpha_0+2}\frac{\alpha}{\alpha_0}c^\alpha B\left(\frac{\alpha_0+2}{2},\frac{\alpha_0+2}{2}\right).$$
We refer to Figure \ref{LSN1:Boxplot} for estimation results comparable with Figure \ref{SN1:Boxplot} with only one line extracted.

\medskip

\begin{figure}[h!]
\begin{center}
\begin{tabular}{ccc}
\hspace{-0.5cm}\includegraphics[width=5.5cm]{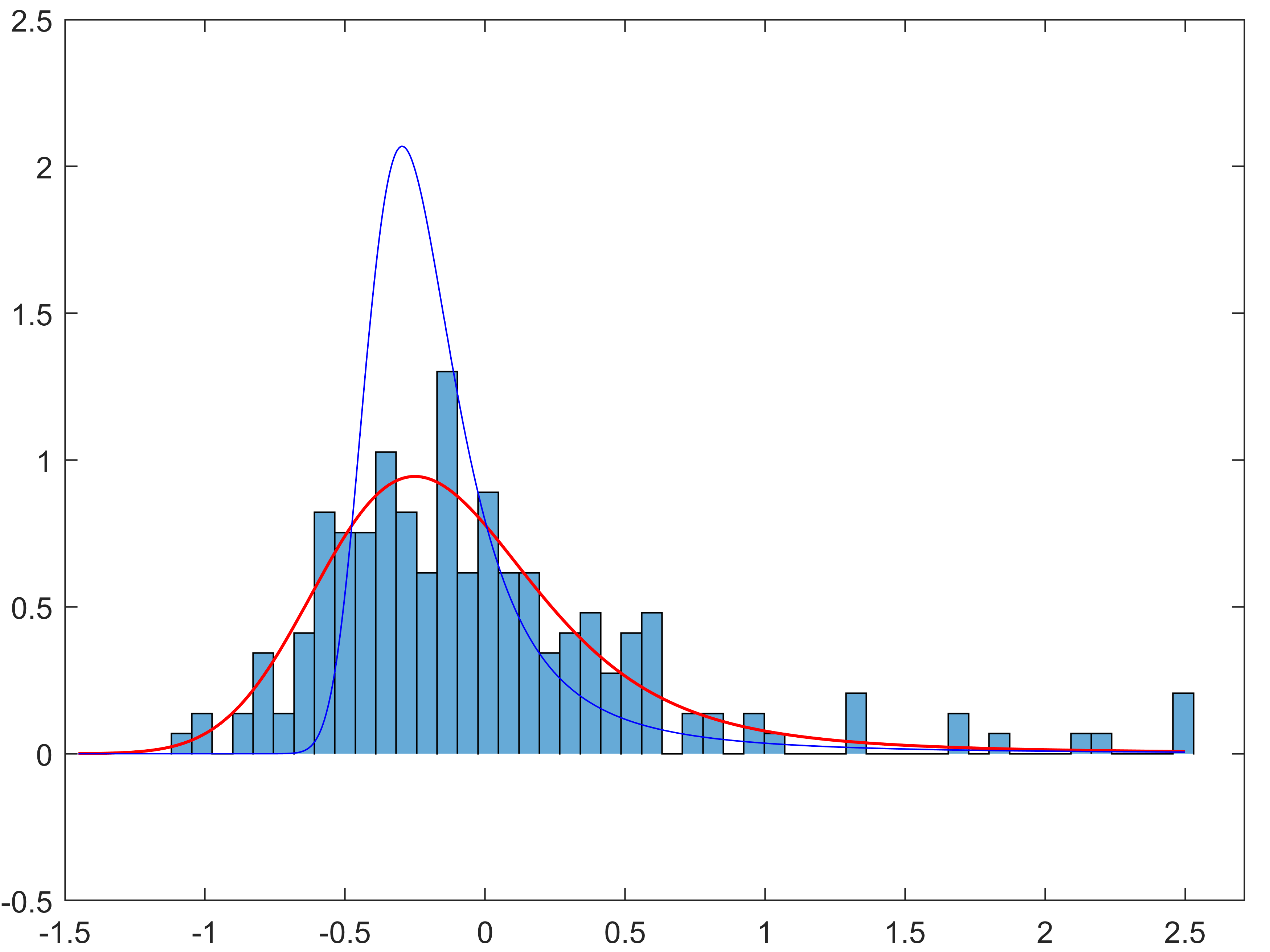} &
\includegraphics[width=5.5cm]{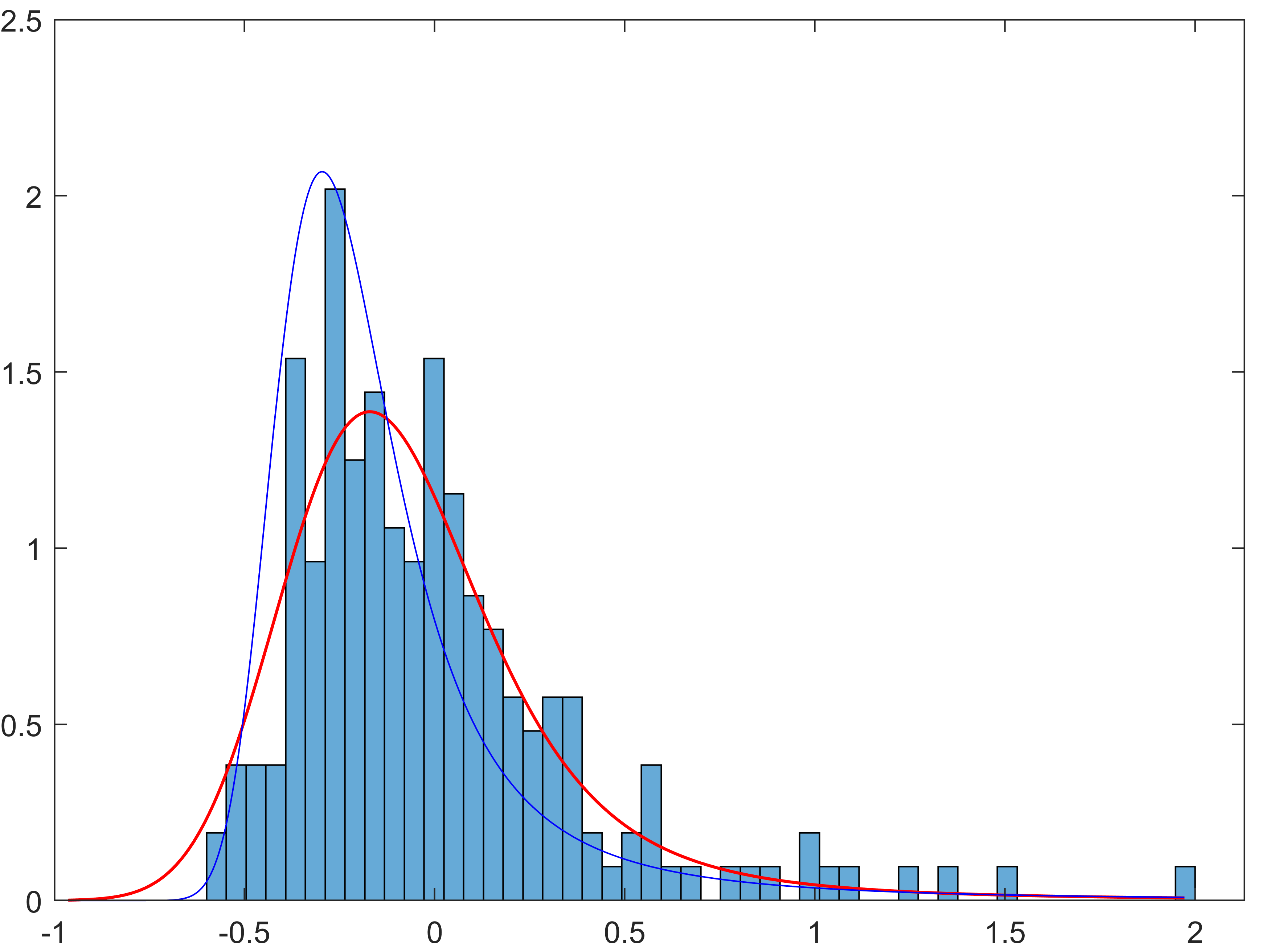} &
\includegraphics[width=5.5cm]{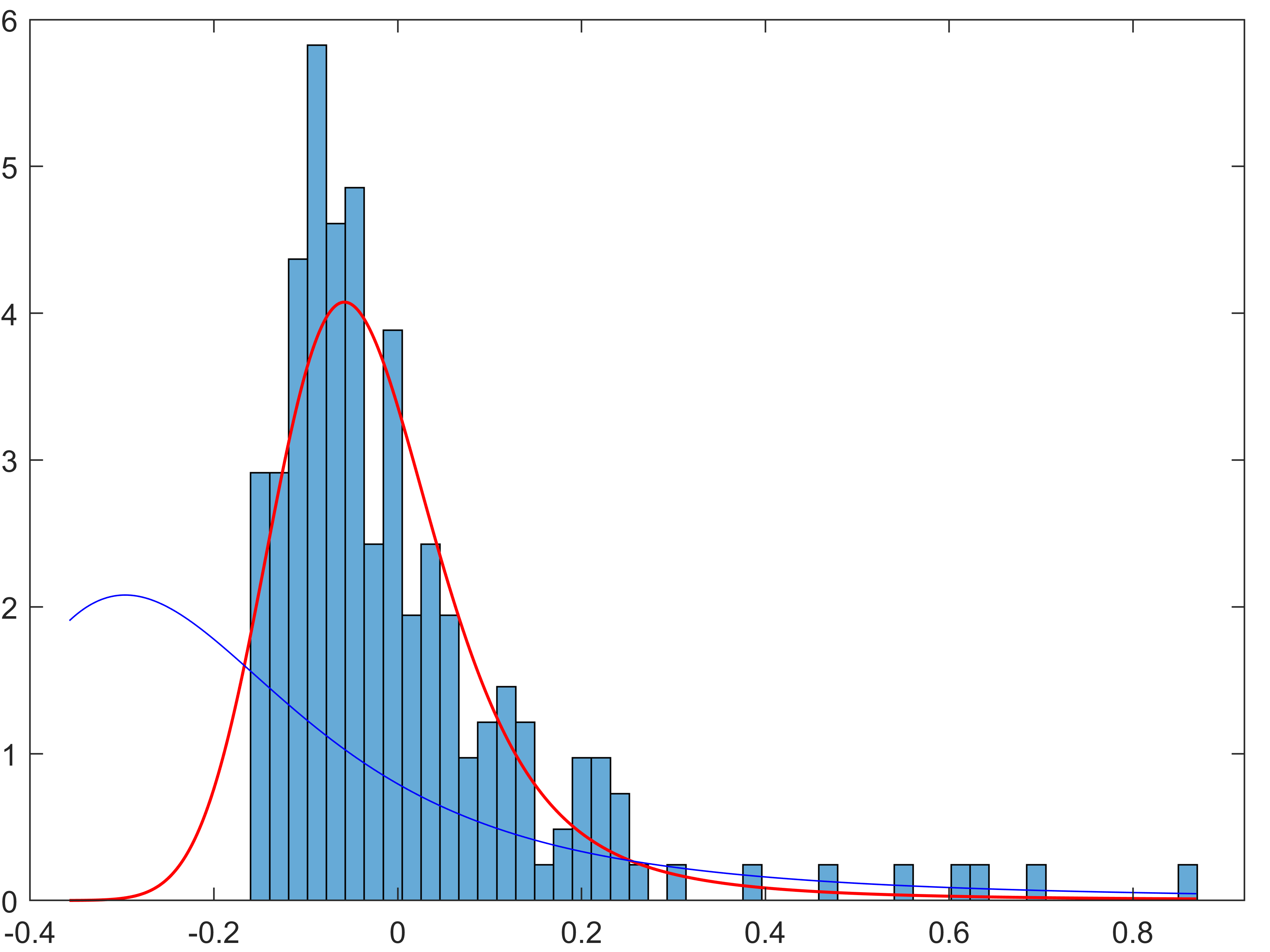}\\
\hspace{-0.5cm}\includegraphics[width=5.5cm]{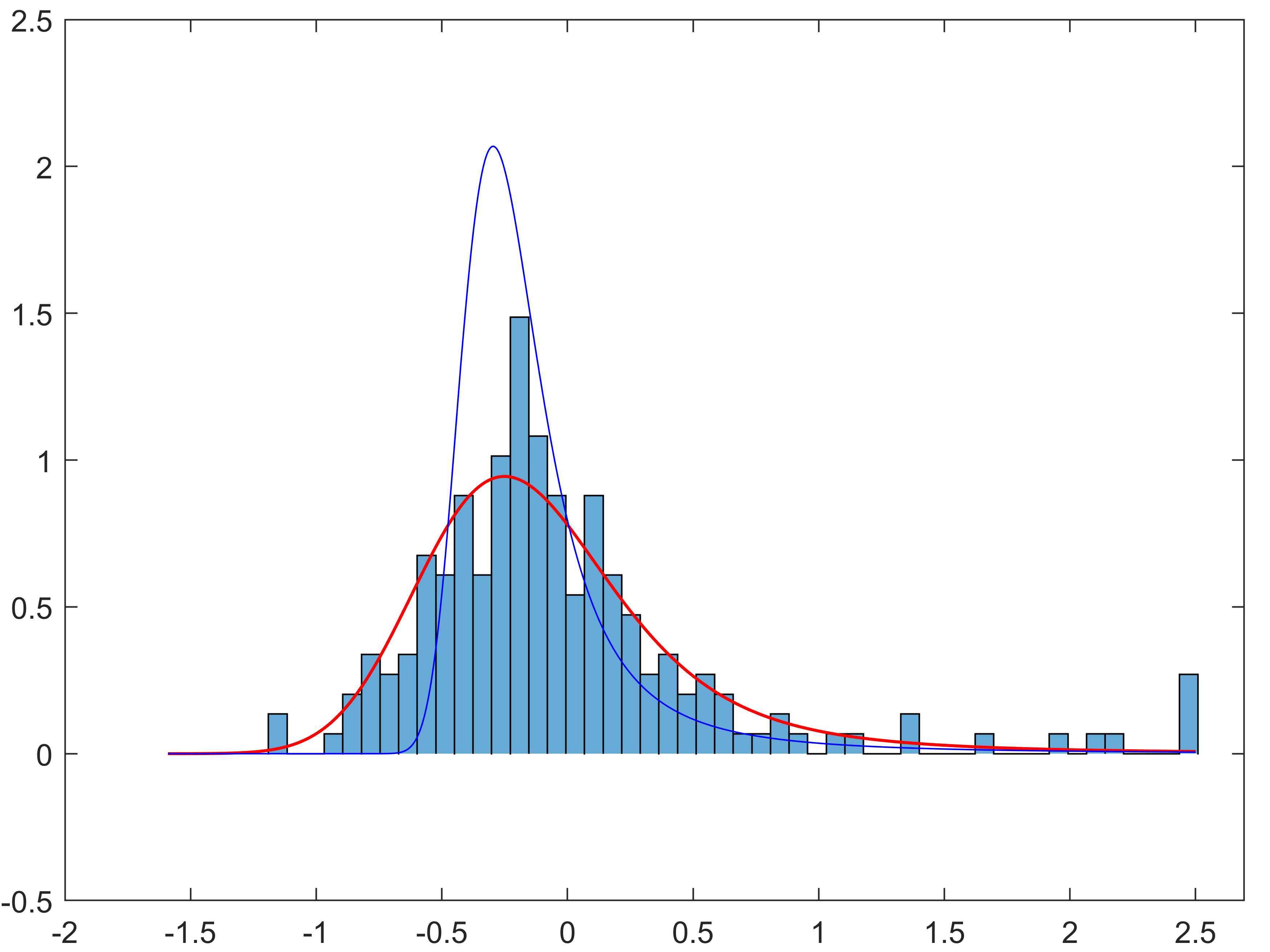} &
\includegraphics[width=5.5cm]{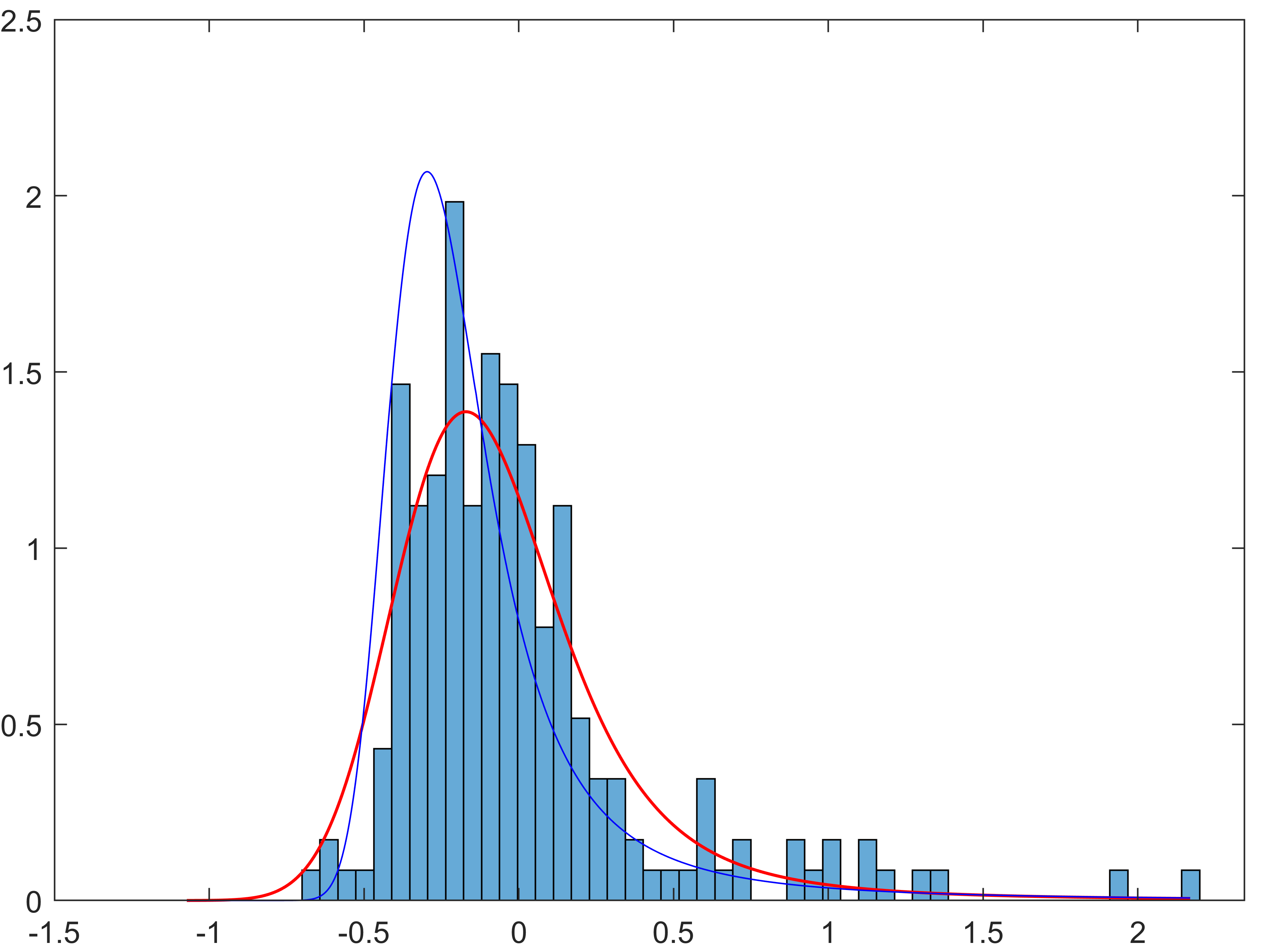} &
\includegraphics[width=5.5cm]{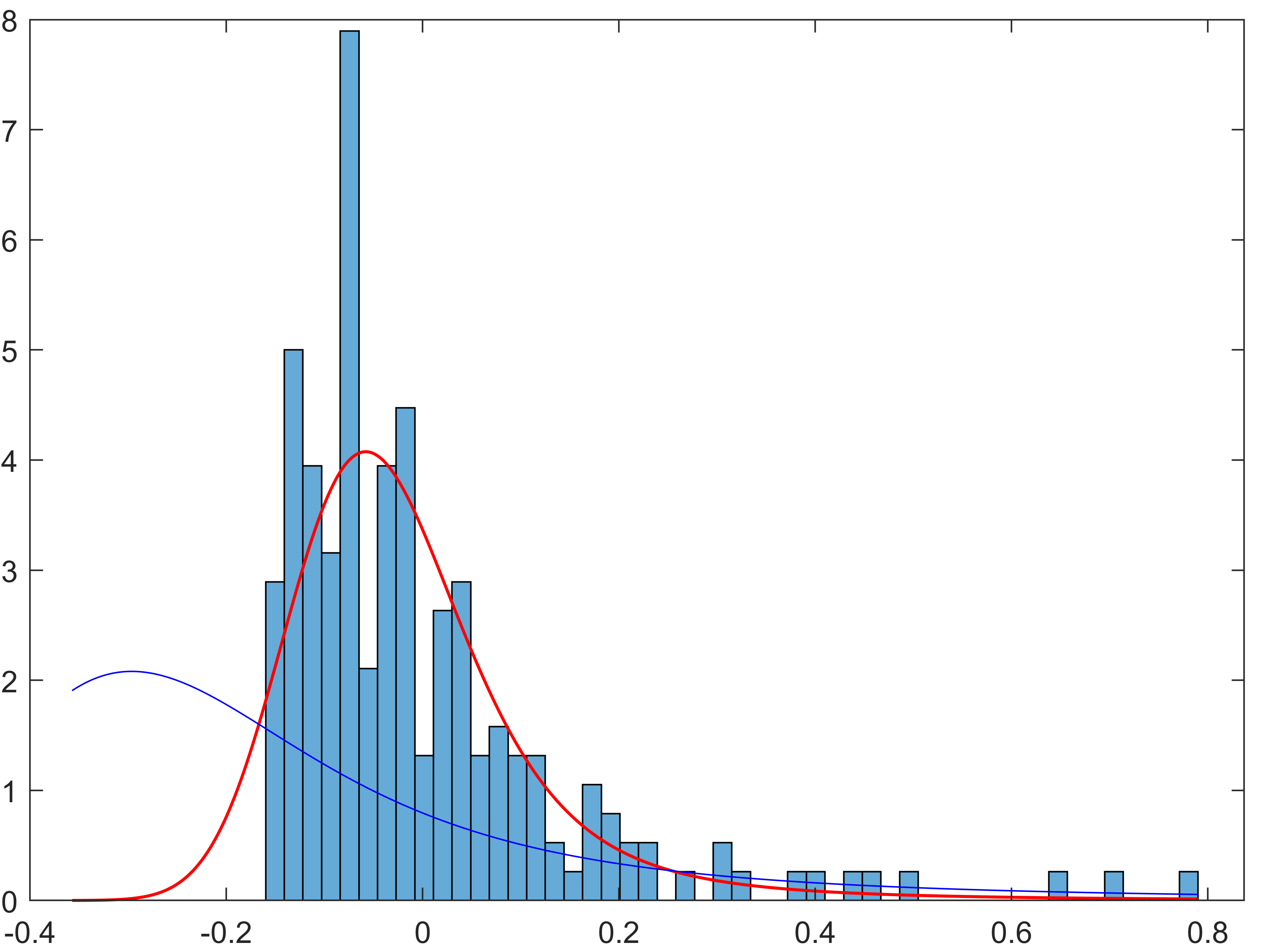}\\
$k=1$&$k=2$&$k=3$
\end{tabular}
\end{center}
\caption{\small Random grain model for $\alpha=1.3$, $a_p=0.005$ and $c=10a_p$. Histograms of the statistic $\lambda^{\nu_0-\nu_0/\alpha_0}\lambda^{\nu-\nu/\alpha} (\widehat Lp_{\lambda,k}-p_k)$ for $k\in\{1,2,3\}$, based on  $200$ independent realizations. The red curves show the theoretical asymptotic  probability distributions in view of \eqref{scale:par0}. The  blue curves show the stable  probability distributions given by \eqref{volumefe1}.}
\label{LSN1:Histo}
\end{figure}

In Figure \ref{LSN1:Histo} we check the asymptotic behavior and compare with the  stable distribution limit obtained for the image.


\newpage

\noindent{\bf Acknowledgments:}
This work has been partly supported by the project ANR MISTIC (ANR-19-CE40-0005).

\end{document}